\documentclass[reqno,12pt,a4paper]{amsart}

\usepackage{mathrsfs}
\usepackage{amssymb}
\usepackage{amsfonts}
\usepackage{latexsym}
\usepackage{amsthm}
\usepackage{graphicx}

\def\lb{\label}

\newcommand{\er}[1]{\textrm{(\ref{#1})}}

\begin{document}


\renewcommand{\theequation}{\arabic{section}.\arabic{equation}}
\theoremstyle{plain}
\newtheorem{theorem}{\bf Theorem}[section]
\newtheorem{lemma}[theorem]{\bf Lemma}
\newtheorem{corollary}[theorem]{\bf Corollary}
\newtheorem{proposition}[theorem]{\bf Proposition}
\newtheorem{definition}[theorem]{\bf Definition}

\newtheorem{remark}[theorem]{\bf Remark}

\def\a{\alpha}  \def\cA{{\mathcal A}}     \def\bA{{\bf A}}  \def\mA{{\mathscr A}}
\def\b{\beta}   \def\cB{{\mathcal B}}     \def\bB{{\bf B}}  \def\mB{{\mathscr B}}
\def\g{\gamma}  \def\cC{{\mathcal C}}     \def\bC{{\bf C}}  \def\mC{{\mathscr C}}
\def\G{\Gamma}  \def\cD{{\mathcal D}}     \def\bD{{\bf D}}  \def\mD{{\mathscr D}}
\def\d{\delta}  \def\cE{{\mathcal E}}     \def\bE{{\bf E}}  \def\mE{{\mathscr E}}
\def\D{\Delta}  \def\cF{{\mathcal F}}     \def\bF{{\bf F}}  \def\mF{{\mathscr F}}
\def\c{\chi}    \def\cG{{\mathcal G}}     \def\bG{{\bf G}}  \def\mG{{\mathscr G}}
\def\z{\zeta}   \def\cH{{\mathcal H}}     \def\bH{{\bf H}}  \def\mH{{\mathscr H}}
\def\e{\eta}    \def\cI{{\mathcal I}}     \def\bI{{\bf I}}  \def\mI{{\mathscr I}}
\def\p{\psi}    \def\cJ{{\mathcal J}}     \def\bJ{{\bf J}}  \def\mJ{{\mathscr J}}
\def\vT{\Theta} \def\cK{{\mathcal K}}     \def\bK{{\bf K}}  \def\mK{{\mathscr K}}
\def\k{\kappa}  \def\cL{{\mathcal L}}     \def\bL{{\bf L}}  \def\mL{{\mathscr L}}
\def\l{\lambda} \def\cM{{\mathcal M}}     \def\bM{{\bf M}}  \def\mM{{\mathscr M}}
\def\L{\Lambda} \def\cN{{\mathcal N}}     \def\bN{{\bf N}}  \def\mN{{\mathscr N}}
\def\m{\mu}     \def\cO{{\mathcal O}}     \def\bO{{\bf O}}  \def\mO{{\mathscr O}}
\def\n{\nu}     \def\cP{{\mathcal P}}     \def\bP{{\bf P}}  \def\mP{{\mathscr P}}
\def\r{\rho}    \def\cQ{{\mathcal Q}}     \def\bQ{{\bf Q}}  \def\mQ{{\mathscr Q}}
\def\s{\sigma}  \def\cR{{\mathcal R}}     \def\bR{{\bf R}}  \def\mR{{\mathscr R}}
\def\S{\Sigma}  \def\cS{{\mathcal S}}     \def\bS{{\bf S}}  \def\mS{{\mathscr S}}
\def\t{\tau}    \def\cT{{\mathcal T}}     \def\bT{{\bf T}}  \def\mT{{\mathscr T}}
\def\f{\phi}    \def\cU{{\mathcal U}}     \def\bU{{\bf U}}  \def\mU{{\mathscr U}}
\def\F{\Phi}    \def\cV{{\mathcal V}}     \def\bV{{\bf V}}  \def\mV{{\mathscr V}}
\def\P{\Psi}    \def\cW{{\mathcal W}}     \def\bW{{\bf W}}  \def\mW{{\mathscr W}}
\def\o{\omega}  \def\cX{{\mathcal X}}     \def\bX{{\bf X}}  \def\mX{{\mathscr X}}
\def\x{\xi}     \def\cY{{\mathcal Y}}     \def\bY{{\bf Y}}  \def\mY{{\mathscr Y}}
\def\X{\Xi}     \def\cZ{{\mathcal Z}}     \def\bZ{{\bf Z}}  \def\mZ{{\mathscr Z}}
\def\be{{\bf e}}
\def\bv{{\bf v}} \def\bu{{\bf u}}
\def\Om{\Omega}
\def\bbD{\pmb \Delta}
\def\mm{\mathrm m}
\def\mn{\mathrm n}

\newcommand{\mc}{\mathscr {c}}

\newcommand{\gA}{\mathfrak{A}}          \newcommand{\ga}{\mathfrak{a}}
\newcommand{\gB}{\mathfrak{B}}          \newcommand{\gb}{\mathfrak{b}}
\newcommand{\gC}{\mathfrak{C}}          \newcommand{\gc}{\mathfrak{c}}
\newcommand{\gD}{\mathfrak{D}}          \newcommand{\gd}{\mathfrak{d}}
\newcommand{\gE}{\mathfrak{E}}
\newcommand{\gF}{\mathfrak{F}}           \newcommand{\gf}{\mathfrak{f}}
\newcommand{\gG}{\mathfrak{G}}           
\newcommand{\gH}{\mathfrak{H}}           \newcommand{\gh}{\mathfrak{h}}
\newcommand{\gI}{\mathfrak{I}}           \newcommand{\gi}{\mathfrak{i}}
\newcommand{\gJ}{\mathfrak{J}}           \newcommand{\gj}{\mathfrak{j}}
\newcommand{\gK}{\mathfrak{K}}            \newcommand{\gk}{\mathfrak{k}}
\newcommand{\gL}{\mathfrak{L}}            \newcommand{\gl}{\mathfrak{l}}
\newcommand{\gM}{\mathfrak{M}}            \newcommand{\gm}{\mathfrak{m}}
\newcommand{\gN}{\mathfrak{N}}            \newcommand{\gn}{\mathfrak{n}}
\newcommand{\gO}{\mathfrak{O}}
\newcommand{\gP}{\mathfrak{P}}             \newcommand{\gp}{\mathfrak{p}}
\newcommand{\gQ}{\mathfrak{Q}}             \newcommand{\gq}{\mathfrak{q}}
\newcommand{\gR}{\mathfrak{R}}             \newcommand{\gr}{\mathfrak{r}}
\newcommand{\gS}{\mathfrak{S}}              \newcommand{\gs}{\mathfrak{s}}
\newcommand{\gT}{\mathfrak{T}}             \newcommand{\gt}{\mathfrak{t}}
\newcommand{\gU}{\mathfrak{U}}             \newcommand{\gu}{\mathfrak{u}}
\newcommand{\gV}{\mathfrak{V}}             \newcommand{\gv}{\mathfrak{v}}
\newcommand{\gW}{\mathfrak{W}}             \newcommand{\gw}{\mathfrak{w}}
\newcommand{\gX}{\mathfrak{X}}               \newcommand{\gx}{\mathfrak{x}}
\newcommand{\gY}{\mathfrak{Y}}              \newcommand{\gy}{\mathfrak{y}}
\newcommand{\gZ}{\mathfrak{Z}}             \newcommand{\gz}{\mathfrak{z}}

\def\ve{\varepsilon}   \def\vt{\vartheta}    \def\vp{\varphi}    \def\vk{\varkappa}

\def\A{{\mathbb A}} \def\B{{\mathbb B}} \def\C{{\mathbb C}}
\def\dD{{\mathbb D}} \def\E{{\mathbb E}} \def\dF{{\mathbb F}} \def\dG{{\mathbb G}} \def\H{{\mathbb H}}\def\I{{\mathbb I}} \def\J{{\mathbb J}} \def\K{{\mathbb K}} \def\dL{{\mathbb L}}\def\M{{\mathbb M}} \def\N{{\mathbb N}} \def\O{{\mathbb O}} \def\dP{{\mathbb P}} \def\R{{\mathbb R}}\def\S{{\mathbb S}} \def\T{{\mathbb T}} \def\U{{\mathbb U}} \def\V{{\mathbb V}}\def\W{{\mathbb W}} \def\X{{\mathbb X}} \def\Y{{\mathbb Y}} \def\Z{{\mathbb Z}}


\def\la{\leftarrow}              \def\ra{\rightarrow}            \def\Ra{\Rightarrow}
\def\ua{\uparrow}                \def\da{\downarrow}
\def\lra{\leftrightarrow}        \def\Lra{\Leftrightarrow}


\def\lt{\biggl}                  \def\rt{\biggr}
\def\ol{\overline}               \def\wt{\widetilde}
\def\ul{\underline}
\def\no{\noindent}


\let\ge\geqslant                 \let\le\leqslant
\def\lan{\langle}                \def\ran{\rangle}
\def\/{\over}                    \def\iy{\infty}
\def\sm{\setminus}               \def\es{\emptyset}
\def\ss{\subset}                 \def\ts{\times}
\def\pa{\partial}                \def\os{\oplus}
\def\om{\ominus}                 \def\ev{\equiv}
\def\iint{\int\!\!\!\int}        \def\iintt{\mathop{\int\!\!\int\!\!\dots\!\!\int}\limits}
\def\el2{\ell^{\,2}}             \def\1{1\!\!1}
\def\sh{\sharp}
\def\wh{\widehat}
\def\bs{\backslash}
\def\intl{\int\limits}

\def\na{\mathop{\mathrm{\nabla}}\nolimits}
\def\sh{\mathop{\mathrm{sh}}\nolimits}
\def\ch{\mathop{\mathrm{ch}}\nolimits}
\def\where{\mathop{\mathrm{where}}\nolimits}
\def\all{\mathop{\mathrm{all}}\nolimits}
\def\as{\mathop{\mathrm{as}}\nolimits}
\def\Area{\mathop{\mathrm{Area}}\nolimits}
\def\arg{\mathop{\mathrm{arg}}\nolimits}
\def\const{\mathop{\mathrm{const}}\nolimits}
\def\det{\mathop{\mathrm{det}}\nolimits}
\def\diag{\mathop{\mathrm{diag}}\nolimits}
\def\diam{\mathop{\mathrm{diam}}\nolimits}
\def\dim{\mathop{\mathrm{dim}}\nolimits}
\def\dist{\mathop{\mathrm{dist}}\nolimits}
\def\Im{\mathop{\mathrm{Im}}\nolimits}
\def\Iso{\mathop{\mathrm{Iso}}\nolimits}
\def\Ker{\mathop{\mathrm{Ker}}\nolimits}
\def\Lip{\mathop{\mathrm{Lip}}\nolimits}
\def\rank{\mathop{\mathrm{rank}}\limits}
\def\Ran{\mathop{\mathrm{Ran}}\nolimits}
\def\Re{\mathop{\mathrm{Re}}\nolimits}
\def\Res{\mathop{\mathrm{Res}}\nolimits}
\def\res{\mathop{\mathrm{res}}\limits}
\def\sign{\mathop{\mathrm{sign}}\nolimits}
\def\span{\mathop{\mathrm{span}}\nolimits}
\def\supp{\mathop{\mathrm{supp}}\nolimits}
\def\Tr{\mathop{\mathrm{Tr}}\nolimits}
\def\BBox{\hspace{1mm}\vrule height6pt width5.5pt depth0pt \hspace{6pt}}


\newcommand\nh[2]{\widehat{#1}\vphantom{#1}^{(#2)}}
\def\dia{\diamond}

\def\Oplus{\bigoplus\nolimits}



\def\qqq{\qquad}
\def\qq{\quad}
\let\ge\geqslant
\let\le\leqslant
\let\geq\geqslant
\let\leq\leqslant
\newcommand{\ca}{\begin{cases}}
\newcommand{\ac}{\end{cases}}
\newcommand{\ma}{\begin{pmatrix}}
\newcommand{\am}{\end{pmatrix}}
\renewcommand{\[}{\begin{equation}}
\renewcommand{\]}{\end{equation}}
\def\eq{\begin{equation}}
\def\qe{\end{equation}}
\def\[{\begin{equation}}
\def\bu{\bullet}

\title[{Magnetic Schr\"odinger operators on periodic discrete graphs}]
{Magnetic Schr\"odinger operators on periodic discrete graphs}

\date{\today}
\author[Evgeny Korotyaev]{Evgeny Korotyaev}
\address{Saint-Petersburg State University, Universitetskaya nab. 7/9, St. Petersburg, 199034, Russia,
\ korotyaev@gmail.com, \
e.korotyaev@spbu.ru,}
\author[Natalia Saburova]{Natalia Saburova}
\address{Northern (Arctic) Federal University, Severnaya Dvina Emb. 17, Arkhangelsk, 163002, Russia,
 \ n.saburova@gmail.com, \ n.saburova@narfu.ru}


\begin{abstract}
We consider magnetic Schr\"odinger operators with periodic magnetic
and electric potentials on periodic  discrete  graphs. The spectrum
of the operators consists of an absolutely continuous part (a~union
of a finite number of non-degenerate bands) plus a finite number of
flat bands, i.e., eigenvalues  of infinite multiplicity. We estimate
the Lebesgue measure of the spectrum in terms of the Betti numbers
and show that these estimates become identities for specific graphs.
We estimate a variation of the spectrum of the Schr\"odinger
operators under a perturbation by a magnetic field in terms of
magnetic fluxes. The proof is based on Floquet theory and a precise
representation of fiber magnetic Schr\"odinger operators constructed
in the paper.

\end{abstract}

\maketitle

\begin{quotation}

\begin{center}
{\bf Table of Contents}
\end{center}

\vskip 6pt

{\footnotesize

1. Introduction \hfill \pageref{Sec1}\ \ \ \ \

2. Main results  \hfill \pageref{Sec2}\ \ \ \ \

3. Direct integrals for magnetic Schr\"odinger operators \hfill
\pageref{Sec3}\ \ \ \ \

4. Proof of the main results \hfill
\pageref{Sec4}\ \ \ \ \

5. Properties of fiber operators and an example \hfill
\pageref{Sec5}\ \ \ \ \

6. Generalized magnetic Schr\"odinger operators \hfill
\pageref{Sec6}\ \ \ \ \

7. Acknowledgments \hfill \pageref{Sec7}\ \ \ \ \

8. References \hfill \pageref{Sec8}\ \ \ \ \
 }
\end{quotation}

\vskip 0.25cm

\section {\lb{Sec1}Introduction}
\setcounter{equation}{0}

\subsection{Introduction.} We discuss spectral properties
of Schr\"odinger operators with periodic magnetic and electric
potentials  on $\Z^d$-periodic discrete graphs, $d\ge 2$, and in
particular, magnetic Laplacians. The spectrum of these operators
consists of an absolutely continuous part (a~union of a finite
number of non-degenerate bands) plus a finite number of flat bands,
i.e., eigenvalues  of infinite multiplicity. There are a lot of
results about such problems, see e.g., \cite{H55}, \cite{Ho76},
\cite{HS89}, \cite{HS01}, \cite{LL93} and the references therein.

A discrete analogue of the magnetic Laplacian on $\R^2$
was originally introduced by Harper \cite{H55}. This discrete
magnetic Laplacian  $\D_\a$ acts on functions $f\in\ell^2(\Z^2), \
n=(n_1,n_2)\in\Z^2$, and is given by:
\[
\lb{haop}
(\D_\a
f)(n)=4f(n)-e^{-iB{n_2\/2}}f(n+e_1)-e^{iB{n_2\/2}}f(n-e_1)-
e^{-iB{n_1\/2}}f(n+e_2)-e^{iB{n_1\/2}}f(n-e_2),
\]
where  $e_1=(1,0), e_2=(0,1)\in \R^2$. The operator $\D_\a$
describes the behavior of an electron moving on the square lattice
$\Z^2$ exposed to a uniform magnetic field in
the so-called tight-binding model \cite{Az64}. The magnetic field
$\cB=B(0,0,1)\in\R^3$ with
amplitude $B\in\R$ is perpendicular to the lattice. The corresponding vector potential $\a$ of the
uniform magnetic field $\cB$ is given by
\[\lb{abe}
\a(\be)=\left\{
\begin{array}{cl}
  -\,{Bn_2\/2}\,, & \textrm{if } \be=(n,n+e_1) \\[6pt]
  -\,{Bn_1\/2}\,, & \textrm{if } \be=(n,n+e_2)
\end{array}\right..
\]
The value $B$ is the magnetic  flux through the unit cell of the
lattice for the magnetic field $\cB$. Note that the discrete
magnetic Laplacian  $\D_\a$ is reduced to the Harper operator in the
discrete Hilbert space $\ell^2(\Z)$. Seemingly it is a very simple
operator but, compared with the magnetic  Laplacian on $\R^2$, its
spectrum is very sensitive to the parameter   $B$ (see \cite{AJ09},
\cite{BS82}, \cite{CEY90}, \cite{Ho76} and the references therein):

  1) if ${B\/2\pi}$ is a rational number, then the spectrum
   $\s(\D_\a)$ of the magnetic Laplacian $\D_\a$ has a \emph{band structure},
   i.e., $\s(\D_\a)$ consists of a finite number of closed intervals;

  2) if ${B\/2\pi}$ is an irrational number, then $\s(\D_\a)$ is a
    \emph{Cantor set} and the graphical presentation of the
    dependence of the spectrum on $B$ shows a
    fractal behavior known as the Hofstadter butterfly.

In a series of papers \cite{HS88}, \cite{HS89}, \cite{HS90} Helffer
and Sj\"ostrand obtained important results in the mathematical
analysis of the magnetic Laplacian $\D_\a$. An algebraic approach to
the operator $\D_\a$ was put forward by Bellissard (for more details
see \cite{Be92}, \cite{Be94}). Note that there are results about the
Hofstadter-type spectrum of the magnetic Laplacians on other planar
graphs (the hexagonal lattice and so on) (see \cite{Hou09},
\cite{Ke92}, \cite{KeR14}
and the references therein).

Discrete magnetic Laplacians on graphs  were introduced
 by Lieb-Loss \cite{LL93} and Sunada \cite{S94}.
 Lieb and Loss
\cite{LL93} characterized the bottom of the spectrum of the discrete
magnetic Laplacian for a bipartite planar graph. Sunada \cite{S94}
considered a discrete magnetic Laplacian with a weak invariance
under a group action on periodic graphs  and gave some criteria
under which the spectrum of the operator has a band structure. After
that, discrete magnetic Schr\"odinger operators on finite  and
infinite graphs have been investigated by many authors.
For example, discrete magnetic Schr\"odinger operators on periodic
graphs were also considered in \cite{HS99a}, \cite{HS99b}. Higuchi
and Shirai \cite{HS99a} obtained the relationship between the
spectrum of the discrete magnetic Schr\"odinger operator on a periodic
graph and that on the corresponding fundamental graph. Also they
proved the analyticity of the bottom of the spectrum with respect to
the magnetic flow and computed the second derivative of the bottom
of the spectrum and represented it in terms of geometry of the
graph. Higuchi and  Shirai \cite{HS99b} gave a condition under which
the weak Bloch property for the magnetic Laplacian holds true, that
is, the set of $\ell^\iy$-eigenvalues is contained in the set of
$\ell^2$-spectrum. Also they investigated spectral properties for
some specific $\Z^d$-periodic graphs $\G$ when $d=\#\cE_*-\#V_*+1$,
where $\#\cE_*$ and $\#V_*$ are the numbers of edges and vertices of
a fundamental graph of $\G$, respectively (see definitions in
subsection 1.2).

Higuchi and Shirai \cite{HS01} studied the behaviour of the bottom
of the spectrum as a function of the magnetic flux. Colin de
Verdi\`{e}re, Torki-Hamza and Truc \cite{CTT11} obtained a condition
under which the magnetic Laplacian on an infinite graph is
essentially self-adjoint.

In our paper we consider the magnetic Laplacians and Schr\"odinger
operators with periodic magnetic and electric potentials on periodic
graphs. The periodicity of magnetic vector potentials guarantees a
band structure of the spectrum and the absence of Cantor spectrum.
Note that in the rational case  ${B\/2\pi}={p\/q}$\,, where $p\in\Z$ and
$q\in\N$ are relatively prime, the vector potential $\a$ defined by
\er{abe} can be considered as a periodic one with the periods
$2qe_1$, $2qe_2$.

We describe now our main goals:

 {\it 1) to estimate the Lebesgue measure of the spectrum and the gaps of
  the magnetic Schr\"odinger operators in terms of
  Betti numbers defined by \er{benu} and electric potentials (see Theorem \ref{T1}).

 2) to estimate a variation of the spectrum of
  the Schr\"odinger operators under a perturbation by
   a magnetic field in terms of magnetic fluxes (see Theorem \ref{Temf}).

 3)  to estimate effective masses associated with the ends of each spectral band for magnetic Laplacians in terms of geometric parameters of the
  graphs (see Theorem \ref{T2}).}

We note that for non-magnetic operators similar estimates were obtained for the Lebesgue measure of the spectrum in \cite{KS14} and for effective masses in \cite{KS16}.

The proof of our results is based on Floquet theory and a precise representation of fiber magnetic Schr\"odinger operators constructed in Theorem \ref{TFR2} and Corollary \ref{TCo0}. This representation of the fiber operators is also the original part of the work. In the proof we use variational estimates for the fiber operators.

\subsection{The definition of magnetic Schr\"odinger operators on periodic graphs.}
Let $\G=(V,\cE)$ be a connected infinite graph, possibly  having
loops and multiple edges, where $V$ is the set of its vertices and
$\cE$ is the set of its unoriented edges. Considering each edge in $\cE$ to have two orientations, we introduce the set $\cA$ of all oriented edges.
An edge starting at a vertex
$u$ and ending at a vertex $v$ from $V$ will be denoted as the ordered pair
$(u,v)\in\cA$ and is said to be \emph{incident} to the vertices.
Vertices $u,v\in V$ will be called \emph{adjacent} and denoted by
$u\sim v$, if $(u,v)\in \cA$.
We define the degree ${\vk}_v$ of
the vertex $v\in V$ as the number of all edges from
$\cA$, starting at $v$. A sequence of directed edges $\cC=(\be_1,\be_2,\ldots,\be_n)$ from $\cA$ is called \emph{a cycle} if the terminus of the edge $\be_s$ coincides with the origin of the edge $\be_{s+1}$ for all $s=1,\ldots,n$ ($\be_{n+1}$ is understood as $\be_1$).

Below we consider locally finite
$\Z^d$-periodic graphs $\G$, $d\geq2$, i.e., graphs satisfying the
following conditions:

{\it 1) $\G$ is equipped with an action of the free abelian group $\Z^d$;

2) the degree of each vertex is finite;

3) the quotient graph  $\G_*=\G/{\Z}^d$ is finite.}

We also call the quotient graph
$\G_*=\G/{\Z}^d$ \emph{the fundamental graph} of the periodic graph $\G$.
If $\G$ is embedded into the space $\R^d$,
the fundamental graph $\G_*$ is a graph on the surface $\R^d/\Z^d$. The fundamental graph $\G_*=(V_*,\cE_*)$ has the vertex set $V_*=V/\Z^d$,
the set $\cE_*=\cE/\Z^d$ of unoriented edges and the set $\cA_*=\cA/\Z^d$ of oriented edges.

\medskip

\no \textbf{Remark.} We do not assume the graph to be embedded into
a Euclidean space. But in many applications there exists such a natural
embedding. The tight-binding approximation is commonly used to
describe the electronic properties of real crystalline structures
(see, e.g., \cite{A76}). This is equivalent to modeling the material
as a discrete graph consisting of vertices (points representing
positions of atoms) and edges (representing chemical bonding of
atoms), by ignoring the physical characters of atoms and bonds that
may be different from one another, see \cite{S13}. The model gives
good qualitative results in many cases. In this case a simple
geometric model is a graph $\G$ embedded into $\R^d$ in such a way
that it is invariant with respect to the shifts by integer vectors
$m\in\Z^d$, which produce an action of $\Z^d$.

\

Let $\ell^2(V)$ be the
Hilbert space of all square summable functions $f:V\to \C$, equipped
with the norm
$$
\|f\|^2_{\ell^2(V)}=\sum_{v\in V}|f(v)|^2<\infty.
$$
Let $\ul\be=(v,u)$ be the inverse edge of $\be=(u,v)\in\cA$. We define the space $\mF_1$ of all periodic 1-forms on a periodic graph $\G$ by
\[\lb{mpin}
\mF_1=\{\a: \cA\ra\R\mid \a(\ul\be\,)=-\a(\be),\qq \a(\be+m)=\a(\be) \qq \textrm{for all } (\be,m)\in \cA\ts\Z^d\},
\]
where $\be+m$ denotes the action of $m\in\Z^d$ on $\be\in\cA$.
In physics a 1-form $\a$ is called \emph{a magnetic vector potential on $\G$}. The number $\a(\be)$, $\be=(u,v)$, is the integral of the magnetic vector potential from the point $u$ to the point $v$.

For each 1-form $\a\in\mF_1$ we define
\emph{the discrete combinatorial magnetic Laplacian} $\D_\a$ on
$f\in\ell^2(V)$ by
\[
\lb{DLO}
 \big(\D_\a f\big)(v)=\sum_{\be=(v,u)\in\cA}\big(f(v)-e^{i\a(\be)}f(u)\big), \qqq
 v\in V.
\]
The sum in \er{DLO} is taken over all oriented edges starting at the vertex $v$.

\

$\bu $ If $\a=0$, then $\D_0$ is just the standard discrete
combinatorial Laplacian $\D$:
\[
\big(\D f\big)(v)=\sum_{(v,u)\in\cA}\big(f(v)-f(u)\big), \qqq v\in V.
\]

$\bu $ If the graph $\G=\Z^2$ and the vector potential $\a$ of the
uniform magnetic field $\cB$ is given by \er{abe}, where the number
${B\/ 2\pi}$ is rational, then $\D_\a$ in \er{DLO} is the operator
defined by \er{haop}.

\

It is well known (see \cite{HS99a}, \cite{HS99b}, \cite{HS01})  that
\emph{the magnetic Laplacian $\D_\a$ is a bounded self-adjoint
operator on $\ell^2(V)$ and its spectrum $\s(\D_\a)$ is a closed
subset in $[0,2\vk_+]$, i.e.:}
\[\lb{vkpl}
\begin{aligned}
\s(\D_\a)\ss [0,2\vk_+],\\
{\rm where } \ \vk_+=\sup_{v\in V}\vk_v<\iy.
\end{aligned}
\]

\

We consider \emph{the magnetic Schr\"odinger operator} $H_\a$ acting on the Hilbert space $\ell^2(V)$ and given by
\[
\lb{Sh}
H_\a=\D_\a+Q,
\]
\[
\lb{Pot}
\big(Q f\big)(v)=Q(v)f(v),\qqq \forall v\in V.
\]
Here and below we assume that the potential $Q$ is real valued and satisfies
$$
Q(v+m)=Q(v), \qqq  \forall\, (v,m)\in V\ts\Z^d,
$$
$v+m$ denotes the action of $m\in\Z^d$ on $v\in V$.

\subsection{Edge indices.}
In order to formulate our results we need to define an {\it edge index}, which was introduced in \cite{KS14}. The indices are important to study the spectrum of the Laplacians and Schr\"odinger operators on periodic graphs, since fiber operators are expressed in terms of edge indices of the fundamental graph (see \er{l2.13}).

Let $\nu=\#V_*$, where $\#A$ is the number of
elements of the set $A$.
We fix any $\nu$ vertices of the periodic graph $\G$, which are not $\Z^d$-equivalent to each other and denote this vertex set by $V_0$.
We will call $V_0$ a \emph{fundamental vertex set of $\G$.} The set $V_0$ is not unique and we may choose this set in different ways. But it is natural to choose the fundamental vertex set $V_0$ in the following way. Let $T=(V_T,\cE_T)$ be a subgraph of the periodic graph $\G$ satisfying the following conditions:

1) \emph{$T$ is a tree, i.e., a connected graph without cycles;}

2) {\it $V_T$ consists of $\n$ vertices of $\G$, which are not $\Z^d$-equivalent to each other.}\\
From now on we assume that the fundamental vertex set $V_0$
coincides with the vertex set $V_T$.

\

\no \textbf{Remark.} Note that such a graph $T$ always exists, since the
periodic graph is connected, and $T$ is not unique.

\

For any vertex $v\in V$ the following unique representation holds true:
\[
\lb{Dv} v=v_0+[v], \qquad v_0\in
V_0,\qquad [v]\in\Z^d.
\]
In other words, each vertex $v$ can be obtained from a vertex $v_0\in V_0$ by the shift by a vector $[v]\in \Z^d$. We will call $[v]$ the \emph{coordinates of the vertex $v$ with respect to the fundamental vertex set $V_0$}.
For any
oriented edge $\be=(u,v)\in\cA$ we define the {\bf edge "index"}
$\t(\be)$ as the integer vector given by
\[
\lb{in}
\t(\be)=[v]-[u]\in\Z^d,
\]
where, due to \er{Dv}, we have
$$
u=u_0+[u],\qquad v=v_0+[v], \qquad u_0,v_0\in V_0,\qquad [u],[v]\in\Z^d.
$$
In general, edge indices  depend on the choice of the set $V_0$.

For example, for the graph $\G$ shown in Fig.\ref{ff.0.11} the index of the edge $(v_1,v_3+a_2)$ is equal to $(0,1)$ and the edge $(v_1,v_4)$ has zero index.

\setlength{\unitlength}{1.0mm}
\begin{figure}[h]
\centering
\unitlength 1mm 
\linethickness{0.4pt}
\ifx\plotpoint\undefined\newsavebox{\plotpoint}\fi 

\begin{picture}(60,50)(0,0)

\put(25.5,44){$\scriptstyle A$}

\multiput(10,10)(4,0){10}{\line(1,0){2}}
\multiput(10,30)(4,0){10}{\line(1,0){2}}
\multiput(10,50)(4,0){10}{\line(1,0){2}}

\multiput(10,10)(0,4){10}{\line(0,1){2}}
\multiput(30,10)(0,4){10}{\line(0,1){2}}
\multiput(50,10)(0,4){10}{\line(0,1){2}}

\put(10,10){\vector(1,0){20.00}}
\put(10,10){\vector(0,1){20.00}}

\put(20,8){$\scriptstyle a_1$}
\put(6,20){$\scriptstyle a_2$}

\put(15,15){\line(1,0){10.00}}
\put(15,15.1){\line(1,0){10.00}}
\put(15,15.2){\line(1,0){10.00}}
\put(15,14.9){\line(1,0){10.00}}
\put(15,14.8){\line(1,0){10.00}}


\put(25,25){\line(1,1){10.00}}

\put(15,15){\line(-1,2){5.00}}
\put(15,15.2){\line(-1,2){5.00}}
\put(15,15.4){\line(-1,2){5.00}}
\put(15,14.8){\line(-1,2){5.00}}
\put(15,14.6){\line(-1,2){5.00}}

\put(15,15){\line(1,1){5.00}}
\put(15,15.2){\line(1,1){5.00}}
\put(15,15.4){\line(1,1){5.00}}
\put(15,14.9){\line(1,1){5.00}}
\put(15,14.8){\line(1,1){5.00}}

\put(25,15){\line(1,2){5.00}}

\put(15,15){\line(1,3){10.00}}

\put(25,15){\line(0,1){10.00}}
\put(25.1,15){\line(0,1){10.00}}
\put(25.2,15){\line(0,1){10.00}}
\put(24.9,15){\line(0,1){10.00}}

\put(15,15){\circle*{1}}
\put(10,25){\circle*{1}}
\put(20,20){\circle*{1}}
\put(25,25){\circle*{1}}
\put(25,15){\circle*{1}}
\put(12,13){$\scriptstyle v_1$}
\put(11,25){$\scriptstyle v_2$}\put(31,25){$\scriptstyle v_2+a_1$}
\put(26,13){$\scriptstyle v_4$}\put(35,32){$\scriptstyle v_1+a_1+a_2$}
\put(24,27){$\scriptstyle v_3$}\put(21,47){$\scriptstyle v_3+a_2$}
\put(18.5,22){$\scriptstyle v_5$}


\put(35,15){\circle{0.7}}
\put(30,25){\circle{0.7}}
\put(40,20){\circle{0.7}}
\put(45,25){\circle{0.7}}
\put(45,15){\circle{0.7}}

\put(15,35){\circle{0.7}}
\put(10,45){\circle{0.7}}
\put(20,40){\circle{0.7}}
\put(25,45){\circle{0.7}}
\put(25,35){\circle{0.7}}


\put(35,35){\circle{0.7}}
\put(30,45){\circle{0.7}}
\put(40,40){\circle{0.7}}
\put(45,45){\circle{0.7}}
\put(45,35){\circle{0.7}}
\end{picture}

\vspace{-0.5cm} \caption{ \footnotesize  A graph $\G$
with the fundamental vertex set $\{v_1,\ldots,v_5\}$; only edges of the fundamental graph $\G_*$ are shown; the vectors $a_1,a_2$ produce an action of $\Z^2$;
the edges of the tree $T$ are marked by bold.} \label{ff.0.11}
\end{figure}
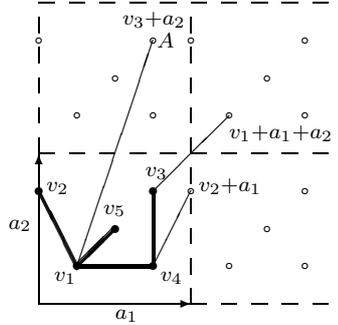

We define two surjections
\[
\gf_V:V\rightarrow V_*=V/\Z^d, \qqq \gf_\cA:\cA\rightarrow\cA_*=\cA/\Z^d,
\]
which map each element to its equivalence class. If $\be$ is an
oriented edge of the graph $\G$, then, by the definition of the
fundamental graph, there is an oriented edge $\be_*=\gf_{\cA}(\be)$
on $\G_*$. For each edge $\be_*\in\cA_\ast$ we define the edge index
$\t(\bf e_*)$ by
\[
\lb{inf}
\t(\bf e_*)=\t(\be).
\]
In other words, edge indices of the fundamental graph $\G_\ast$  are
induced by edge indices of the periodic graph $\G$.
An index of a
fundamental graph edge with respect to the fixed fundamental vertex set $V_0$ is uniquely determined by \er{inf}, since
$$
\t(\be+m)=\t(\be),\qqq \forall\, (\be,m)\in\cA \ts \Z^d.
$$

From the definitions \er{in}, \er{inf} of the
edge index we have
\[\lb{inin}
\t(\ul\be)=-\t(\be) \qq \textrm{ for each } \be\in\cA \textrm{ and for each }\be\in\cA_*;
\]
and \textbf{all edges of the tree $T$ have zero indices},
i.e.,
\[\lb{iesp}
\t(\be)=0, \qqq \forall \ \be\in\cE_{T}.
\]

\subsection{Betti number, spanning trees and magnetic fluxes}
We recall the definitions of the Betti number and spanning trees, which will be used in the formulation of our results.

$\bullet$ The \emph{Betti number} $\b$ of a finite connected graph $\G_*=(V_*,\cE_*)$ is defined as
\[\lb{benu}
\b=\#\cE_*-\# V_*+1.
\]
Note that the Betti number $\b$ can also be defined in one of the following ways:
\begin{itemize}
  \item[i)] as the number of edges that have to be removed from $\cE_*$ (without reducing the number of vertices) to turn $\G_*$ into a tree;
  \item[ii)] as the dimension of the cycle space of the graph $\G_*$
\end{itemize}
(see properties of spanning trees below).

\

$\bullet$  A \emph{spanning tree} $T_*=(V_*,\cE_{T_*})$ of a finite connected graph $\G_*=(V_*,\cE_*)$ is a connected subgraph of $\G_*$ which has no cycles and contains all vertices of $\G_*$.

\

We introduce the set $\cS$ of all edges from $\cE_*$ that do not belong to the spanning tree $T_*$ and equip each edge of $\cS$ with some orientation. We denote by $\ul\cS$ the set of their inverse edges, i.e.,
\[
\lb{ulS}
\cS=\cE_*\sm\cE_{T_*},\qqq \ul\cS=\{\be\in\cA_*\mid \ul\be\in\cS\}.
\]

We recall some \emph{properties of spanning trees} of connected graphs
(see, e.g., Lemma 5.1 and Theorem 5.2 in \cite{B74}):

\emph{1) The set $\cS$ contains $\b$ edges, where $\b$ is the Betti
number defined by \er{benu}.}

\emph{2) For any edge $\be\in\cS$ there exists a unique cycle $\cC_\be$ containing only $\be$ and edges of $T_*$.}

\emph{3) The set of all such cycles $(\cC_\be)_{\be\in\cS}$ forms a basis of the cycle space and the number of independent cycles of the fundamental graph $\G_*$ is $\b$.}

\

\textbf{Remark.} The definitions of the Betti number and spanning trees and their properties hold true for any finite connected graph $\G_*=(V_*,\cE_*)$, which is not necessarily a fundamental graph of some periodic one.

\

For a given magnetic Laplacian $\D_\a$ the magnetic vector potential $\a$ is defined
up to a gauge transformation. Therefore, we define a \emph{magnetic flux}, which is invariant under the gauge transformation.


We recall that $T=(V_T,\cE_T)$ is a connected subgraph of the periodic graph $\G$ with no cycles and with $\n$ vertices which are not $\Z^d$-equivalent to each other, where $\n$ is the number of the fundamental graph vertices. Then the graph $T_*=T/{\Z}^d$ is a spanning tree of
the fundamental graph $\G_*=\G/{\Z}^d$.
Due to the property 2) of spanning trees, for each $\be\in\cS$, where $\cS$ is defined in \er{ulS}, there exists a unique cycle
$\cC_\be$ containing only $\be$ and edges of $T_*$. For the cycle
$\cC_\be$ we define the \emph{magnetic flux} of $\a$ by
\[\lb{mafl'}
\phi_\a(\be)\equiv\phi_\a(\cC_\be)=\Big(\sum_{\wt\be\in\cC_\be}\a(\wt\be)\Big) \; \mathrm{mod} \; 2\pi,\qqq \phi_\a(\be)\in(-\pi,\pi].
\]

\textbf{Example.} For the graph $\G_*$ shown in Fig.\ref{ffS'}\emph{a} we can choose the spanning trees $T_*$ and $\wt T_*$ (Fig.\ref{ffS'} \emph{b,c}).
The set $\cS$ consists of three edges $\be_1,\be_2,\be_3$ (they are shown in Fig.\ref{ffS'} \emph{b,c} by the dotted lines) and depends on the choice of the spanning tree. The Betti number $\b$ defined by \er{benu} is equal to 3 and does not depend on the set $\cS$.

\setlength{\unitlength}{1.0mm}
\begin{figure}[h]
\centering
\unitlength 1.0mm 
\linethickness{0.4pt}
\begin{picture}(100,30)

\put(10,10){\circle{1}}
\put(10,30){\circle{1}}
\put(10,20){\circle{1}}
\put(0,20){\circle{1}}
\put(20,20){\circle{1}}
\put(10,10){\line(1,1){10.00}}
\put(10,10){\line(-1,1){10.00}}
\put(10,30){\line(1,-1){10.00}}
\put(10,30){\line(-1,-1){10.00}}

\put(10,30){\line(0,-1){10.00}}
\put(0,20){\line(1,0){20.00}}

\put(50,10){\circle{1}}
\put(50,30){\circle{1}}
\put(50,20){\circle{1}}
\put(40,20){\circle{1}}
\put(60,20){\circle{1}}
\put(50,10){\line(1,1){10.00}}
\put(50,10.2){\line(1,1){10.00}}
\put(50,9.8){\line(1,1){10.00}}
\put(50,10.1){\line(1,1){10.00}}
\put(50,9.9){\line(1,1){10.00}}

\put(50,30){\line(0,-1){10.00}}
\put(50.1,30){\line(0,-1){10.00}}
\put(49.9,30){\line(0,-1){10.00}}
\put(40,20.1){\line(1,0){20.00}}
\put(40,19.9){\line(1,0){20.00}}
\put(40,20){\line(1,0){20.00}}

\qbezier[20](50,10)(45,15)(40,20)
\qbezier[20](50,30)(55,25)(60,20)
\qbezier[20](50,30)(45,25)(40,20)

\put(42,13){$\be_1$}
\put(41,26){$\be_2$}
\put(55,26){$\be_3$}
\put(90,10){\circle{1}}
\put(90,30){\circle{1}}
\put(90,20){\circle{1}}
\put(80,20){\circle{1}}
\put(100,20){\circle{1}}
\put(90,10){\line(1,1){10.00}}
\put(90,10.2){\line(1,1){10.00}}
\put(90,9.8){\line(1,1){10.00}}
\put(90,10.1){\line(1,1){10.00}}
\put(90,9.9){\line(1,1){10.00}}

\put(90,30){\line(0,-1){10.00}}
\put(90.1,30){\line(0,-1){10.00}}
\put(89.9,30){\line(0,-1){10.00}}
\put(90,20.1){\line(1,0){10.00}}
\put(90,19.9){\line(1,0){10.00}}
\put(90,20){\line(1,0){10.00}}

\qbezier[15](80,20)(85,20)(90,20)
\qbezier[20](90,10)(85,15)(80,20)
\qbezier[20](90,30)(95,25)(100,20)

\put(90,30.2){\line(-1,-1){10.00}}
\put(90,29.8){\line(-1,-1){10.00}}
\put(90,30.1){\line(-1,-1){10.00}}
\put(90,29.9){\line(-1,-1){10.00}}
\put(90,30){\line(-1,-1){10.00}}
\put(82,13){$\be_1$}
\put(84,21){$\be_2$}
\put(95,26){$\be_3$}

\put(-3,10){(\emph{a})}
\put(33,10){(\emph{b})}
\put(73,10){(\emph{c})}

\put(0,26){$\G_*$}
\put(33,26){$T_*$}
\put(80,26){$\wt T_*$}
\end{picture}

\vspace{-10mm}
\caption{\footnotesize  \emph{a}) A fundamental graph $\G_*$;\quad \emph{b}),\emph{c}) the spanning trees $T_*$ and $\wt T_*$, $\cS=\{\be_1,\be_2,\be_3\}$, $\b=3$.} \label{ffS'}
\end{figure}

\section{\lb{Sec2}Main results}
\setcounter{equation}{0}

\subsection{Floquet decomposition of Schr\"odinger operators}
We introduce the Hilbert space
\[\lb{Hisp}
\mH=L^2\Big(\T^d,{d\vt\/(2\pi)^d}\,,\cH\Big)=\int_{\T^d}^{\os}\cH\,{d\vt
\/(2\pi)^d}\,, \qq \cH=\ell^2(V_*), \qq \T^d=\R^d/(2\pi\Z)^d,
\]
i.e., a constant fiber direct integral equipped with the norm
$$
\|g\|^2_{\mH}=\int_{\T^d}\|g(\vt,\cdot)\|_{\ell^2(V_*)}^2\frac{d\vt
}{(2\pi)^d}\,,
$$
where the function $g(\vt,\cdot)\in\cH$ for almost all $\vt\in\T^d$.

\begin{theorem}[\textbf{Magnetic fluxes representation}]
\label{TFR2}
For each 1-form $\a\in\mF_1$ the magnetic Schr\"odinger operator $H_\a=\D_\a+Q$ on $\ell^2(V)$ has the following decomposition into a constant fiber direct integral
\[
\lb{raz1}
\begin{aligned}
& \ell^2(V)={1\/(2\pi)^d}\int^\oplus_{\T^d}\ell^2(V_*)\,d\vt ,\qqq
\mU H_\a \mU^{-1}={1\/(2\pi)^d}\int^\oplus_{\T^d}H_\a(\vt)d\vt,
\end{aligned}
\]
where the unitary operator $\mU:\ell^2(V)\to\mH$ is a
composition of the Gelfand type transformation and a gauge
transformation (see the precise formulas \er{5001} and \er{Uvt}). Here
the fiber magnetic Schr\"odinger operator $H_\a(\vt)$  and  the
fiber magnetic Laplacian $\D_\a(\vt)$ are given by
\[
\label{Hvt}
H_\a(\vt)=\D_\a(\vt)+Q,\qqq \forall\,\vt\in \T^d,
\]
\[
\label{l2.13}
\begin{aligned}
\big(\D_\a(\vt)f\big)(v)=\vk_vf(v)
-\sum_{\be=(v,u)\in\cA_*}e^{i(\a_*(\be)+\lan\t(\be),\,\vt\ran)}f(u), \qqq v\in V_*,
\end{aligned}
\]
where the modified 1-form $\a_*\in\mF_1$ is uniquely defined by
\[\lb{cat}
\a_*(\be)=\left\{
\begin{array}{ll}
 \phi_\a(\be), &  \textrm{ if } \, \be\in \cS\\[6pt]
  0, \qq & \textrm{ if } \, \be\notin\cS\cup\ul\cS\\
\end{array}\right.,
\]
the magnetic flux $\phi_\a(\be)$ is given by \er{mafl'};
$\t(\be)$ is the index of the edge $\be$ defined by \er{in}, \er{inf}; $\cS$ and $\ul\cS$ are
defined by \er{ulS}, and $\lan\cdot\,,\cdot\ran$ denotes the
standard inner product in $\R^d$.
\end{theorem}

 {\bf Remarks.} 1) The  modified  magnetic vector potential
$\a_*\in\mF_1$ on each edge $\be\in\cS$ coincides with the flux
$\phi_\a(\be)$ through the cycle $\cC_\be$.

2) Note that the decomposition of the discrete magnetic
Schr\"odinger operators on periodic graphs into the constant fiber
direct integral \er{raz1} (without an exact form of fiber operators)
was discussed by Higuchi and Shirai \cite{HS99a}. The precise form
of the fiber Laplacian $\D_\a(\vt)$ defined by \er{l2.13} is
important to study spectral properties of the magnetic Laplacians
and Schr\"odinger operators acting on periodic graphs (see the proof of
Theorems \ref{T1} -- \ref{T2}). The precise forms of the fiber Laplacian at $\a=0$ and of fiber metric Laplacians on periodic graphs
were determined in \cite{KS14}, \cite{KS15}.

3) From Theorem \ref{TFR2} it follows that two Schr\"odinger
operators with the same potential and the same magnetic flux through
every basic cycle are unitarily equivalent. This property for the
magnetic Schr\"odinger operators on a locally finite graph was
proved in \cite{LL93}, \cite{CTT11}, \cite{HS01}. In particular, if the magnetic flux of $\a$ is zero for any cycle on $\G_*$, then
the magnetic Schr\"odinger operator $H_\a=\D_\a+Q$ is unitarily equivalent to the Schr\"odinger operator $H_0=\D_0+Q$ without a magnetic field.

4) The modified 1-form $\a_*$ given by \er{cat} depends on the choice of the spanning tree $T_*$.

\

In \er{l2.13}, \er{cat} the fiber magnetic Laplacian $\D_\a(\vt)$ depends on $\b$, generally speaking, non-zero independent magnetic fluxes $\big(\phi_\a(\be)\big)_{\be\in\cS}$. Now we show that using a simple change of variables we can reduce the number of these independent parameters to $\b-d$. In particular, if $\b=d$, then the fiber Laplacian does not depend on the magnetic fluxes.

\begin{corollary}[\textbf{Minimal magnetic fluxes representation}]
\label{TCo0} There exist $\vt_0\in\T^d$ and edges
$\be_1,\ldots,\be_d\in\cS$  with linearly independent indices
$\t(\be_1),\ldots,\t(\be_d)$ defined by \er{in}, \er{inf}
such that the fiber Laplacian $\D_\a(\vt)$ given by \er{l2.13} in the new variables $\wt\vt=\vt-\vt_0$ has the form
\[
\label{wtl2.13}
\begin{aligned}
\big(\D_\a(\wt\vt+\vt_0)f\big)(v)=\vk_vf(v)
-\sum_{\be=(v,u)\in\cA_*}e^{i(\wt\a(\be)+\lan\t(\be),\,\wt\vt\,\ran)}f(u), \qqq v\in V_*,
\end{aligned}
\]
where the modified 1-form $\wt\a\in\mF_1$ is defined by
\[\lb{wtcat}
\wt\a(\be)=\left\{
\begin{array}{ll}
 \phi_\a(\be)+\lan\t(\be),\vt_0\ran, &  \textrm{ if } \, \be\in \wt\cS\\[6pt]
  0, \qq & \textrm{ if } \, \be\notin\wt\cS\cup\ul{\wt\cS}\\
\end{array}\right.,
\]
the magnetic flux $\phi_\a(\be)$ is given by \er{mafl'};
\[
\lb{wulS}
\wt\cS=\cS\sm\{\be_1,\ldots,\be_d\},\qqq \ul{\wt\cS}=\{\be\in\cA_*\mid \ul\be\in\wt\cS\}.
\]
In particular, if the Betti
number $\b$ defined by \er{benu} is equal to
$d$, then the magnetic Schr\"odinger operator
$H_\a$ is unitarily equivalent to the Schr\"odinger operator $H_0$
without a magnetic field.
\end{corollary}

\no \textbf{Remarks.} 1) Higuchi and Shirai \cite{HS99b} show that if
$\b=d$, then  the magnetic Schr\"odinger operator $H_\a$ is
unitarily equivalent to the Schr\"odinger operator $H_0$ without a
magnetic field.  Their proof is based on homology theory. Our proof
is based on a simple change of variables.

2) The hexagonal
lattice and the $d$-dimensional lattice with the minimal fundamental graphs are examples of
periodic graphs with $\b=d$. It is known that for any magnetic
vector potential $\a\in\mF_1$ the spectrum of the magnetic Laplacian
$\D_\a$ on these graphs is given by $\s(\D_\a)=[0,2\vk_+]$, where
$\vk_+$ is the degree of each vertex of the graph, i.e., the
spectrum does not depend on the magnetic potential $\a$.
If the fundamental graphs are not minimal, then $\b>d$ and, in
general, the spectrum of the magnetic Laplacian $\D_\a$ depends on $\a$.

\subsection{Spectrum of the magnetic Schr\"odinger operator.}
Theorem \ref{TFR2} and standard arguments (see Theorem XIII.85 in
\cite{RS78}) describe the spectrum of the magnetic Schr\"odinger operator
$H_\a=\D_\a+Q$. Each fiber operator $H_\a(\vt)$, $\vt\in\T^d$,
has $\n$ eigenvalues $\l_{\a,n}(\vt)$, $n\in\N_\n=\{1,\ldots,\n\}$, $\n=\# V_*$, which are labeled
in increasing order (counting multiplicities) by
\[
\label{eq.3} \l_{\a,1}(\vt)\leq\l_{\a,2}(\vt)\leq\ldots\leq\l_{\a,\nu}(\vt),
\qqq \forall\,\vt\in\T^d.
\]
Since $H_\a(\vt)$ is
self-adjoint and analytic in $\vt\in\T^d$, each $\l_{\a,n}(\cdot)$, $n\in\N_\n$, is a real and piecewise analytic function on the torus $\T^d$ and creates the \emph{spectral band} $\s_n(H_\a)$ given by
\[
\lb{ban.1}
\s_n(H_\a):=\s_{\a,n}=[\l_{\a,n}^-,\l_{\a,n}^+]=\l_{\a,n}(\T^d).
\]
Thus, the spectrum of the operator $H_\a$ on the periodic graph $\G$ is
given by
\[\lb{spec}
\s(H_\a)=\bigcup_{\vt\in\T^d}\s\big(H_\a(\vt)\big)=\bigcup_{n=1}^{\nu}\s_n(H_\a).
\]
Note that if $\l_{\a,n}(\cdot)= C_{\a,n}=\const$ on some subset of $\T^d$ of
positive Lebesgue measure, then  the operator $H_\a$ on $\G$ has the
eigenvalue $C_{\a,n}$ of infinite multiplicity. We call $C_{\a,n}$
a \emph{flat band}.

Thus, the spectrum of the magnetic Schr\"odinger operator
$H_\a$ on the periodic graph $\G$ has the form
\[
\lb{r0}
\s(H_\a)=\s_{ac}(H_\a)\cup \s_{fb}(H_\a).
\]
Here $\s_{ac}(H_\a)$ is the absolutely continuous spectrum, which is a
union of non-degenerate intervals, and $\s_{fb}(H_\a)$ is the set of
all flat bands (eigenvalues of infinite multiplicity). An open
interval between two neighboring non-degenerate spectral bands is
called a \emph{spectral gap}.

The eigenvalues of the fiber magnetic Laplacian $\D_\a(\vt)$  will be denoted by
$\l^0_{\a,n}(\vt)$, $n\in\N_\n$. The spectral bands $\s_n(\D_\a)$, $n\in\N_{\n}$, for the magnetic Laplacian $\D_\a$ have the form
\[
\lb{ban0} \s_n(\D_\a)=[\l_{\a,n}^{0-},\l_{\a,n}^{0+}]=\l_{\a,n}^0(\T^d).
\]

\no \textbf{Remark.} From \er{ban.1} it follows that
\[\lb{HiSi}
\l_{\a,1}^{-}\leq\l_{\a,1}(0)
\]
for any $\a\in \mF_1$. Note that if there is no magnetic field, that is $\a=0$, Sy and Sunada \cite{SS92} proved
that $\l_{0,1}^{-}=\l_{0,1}(0)$. However, the equality in \er{HiSi} does not hold for general $\a$, since for some specific graphs we have the strict inequality (see Examples 5.2 -- 5.6 in \cite{HS01}).

\subsection{Estimates of the Lebesgue measure of the spectrum}
Now we estimate the Lebesgue measure of the spectrum of the magnetic Schr\"odinger operator in terms of the Betti number and the Lebesgue measure of the gaps in terms of the Betti number and electric potentials.

\begin{theorem}
\lb{T1} i) The Lebesgue measure $|\s(H_\a)|$ of the spectrum of the magnetic Schr\"odinger operator $H_\a=\D_\a+Q$ satisfies
\[
\lb{eq.7'}
|\s(H_\a)|\le \sum_{n=1}^{\n}|\s_n(H_\a)|\le 4\b,
\]
where $\b$ is the Betti number defined by \er{benu}.
Moreover, if there exist $s$  spectral gaps \linebreak ${\g_1(H_\a),\ldots,\g_s(H_\a)}$ in the spectrum $\s(H_\a)$, then the following estimates hold true:
\[\lb{GEga}
\begin{aligned}
\sum_{n=1}^s|\g_n(H_\a)|\ge
\l^+_{\a,\n}-\l_{\a,1}^--4\b\ge C_0-4\b,\\
C_0=|\l^{0+}_{\a,\n}-\l^{0-}_{\a,1}-q_\bu|, \qqq q_\bu=\max_{v\in V_*} Q(v)-\min_{v\in V_*} Q(v),
\end{aligned}
\]
where $\l^{0+}_{\a,\n}$, $\l^{0-}_{\a,1}$ and $\l^{+}_{\a,\n}$, $\l^{-}_{\a,1}$ are the upper and lower endpoints of the spectrum of the Laplacian $\D_\a$ and the Schr\"odinger operator $H_\a$, respectively.

ii) The estimates \er{eq.7'} and the first estimate in \er{GEga} become identities for some classes of graphs, see \er{sp2}.
\end{theorem}

\no \textbf{Remarks.} 1) There exists a $\Z^d$-periodic graph $\G$, such that the total
length of all spectral bands of the magnetic Schr\"odinger operators
$H_\a=\D_\a+Q$ on the graph $\G$ \textbf{depends on neither the
potential $Q$ nor the magnetic potential $\a$} (see Proposition
\ref{TG1}).

2) The Lebesgue measure $|\s(H_\a)|$ of the spectrum of $H_\a$ (on
specific graphs) can be arbitrary large (see Proposition \ref{TG1}).

\

Now we estimate a variation of the spectrum of
the Schr\"odinger operators under a perturbation by
a magnetic field in terms of magnetic fluxes.

\begin{theorem}
\lb{Temf} Let $\a\in \mF_1$ and let the corresponding magnetic Schr\"odinger operator
$H_\a=\D_\a+Q$ have the spectral bands $\s_{\a,n}$
given by \er{ban.1}.
Then for any $\a^o\in \mF_1$ all corresponding band ends
$\l_{\a^o,n}^\pm$ satisfy
\[
\lb{emf1} \L_1\le \l_{\a^o,n}^\pm-\l_{\a,n}^\pm\le \L_\n,
\]
\[
\lb{emf2} \big||\s_{\a^o,n}|-|\s_{\a,n}|\big|\le\L_\n-\L_1,
\]
where
\[
\lb{emf3'}\L_1=\min_{\vt\in\T^d} \l_1(X_{\a^o,\a}(\vt)),\qqq
\L_\n=\max_{\vt\in\T^d} \l_\n(X_{\a^o,\a}(\vt)),
\]
and $X_{\a^o,\a}(\cdot)=H_{\a^o}(\cdot)-H_{\a}(\cdot)$. Moreover,
$\L_1$ and $\L_\n$ satisfy the following estimates:
\[
\begin{aligned}
\lb{es30} \max\{|\L_1|,|\L_\n|\}\le C_{\a^o,\a},\qqq \L_\n-\L_1\le
2C_{\a^o,\a},
\end{aligned}
\]
where
\[
\begin{aligned}
\lb{es3} C_{\a^o,\a}=2\max_{u\in V_*}\sum_{\be=(u,v)\in\cS\cup \ul\cS}|\sin x_{\be}|, \qqq x_{\be}={1\/2}\big(\phi_{\a^o}(\be)-\phi_\a(\be)\big),
\end{aligned}
\]
$\cS$ and $\ul\cS$ are
defined by \er{ulS}, and the magnetic flux $\phi_\a(\be)$ is given by \er{mafl'}.
\end{theorem}

{\bf Remark.} The magnetic Schr\"odinger operators depend on magnetic potentials, but we obtain the estimates of a variation of the spectrum in terms of the difference of magnetic fluxes
only.

\subsection{Effective masses for magnetic Laplacians.}
Let $\l_\a(\vt)$, $\vt\in\T^d$, be a band function of the magnetic Laplacian
$\D_\a$ and let $\l_\a(\vt)$ have a minimum (maximum) at some point
$\vt_0$. Assume that $\l_\a(\vt_0)$ is a simple eigenvalue of
$\D_\a(\vt_0)$. Then the eigenvalue $\l_\a(\vt)$ has the
Taylor series as $\vt=\vt_0+\ve\o$, $\o=(\o_\a)_{\a=1}^d\in\S^{d-1}$, $\ve=|\vt-\vt_0|\to0$:
\[
\label{lam}
\begin{aligned}
\l_\a(\vt)=\l_\a(\vt_0)+\ve^2\m_\a(\o)+O(\ve^3),\qqq
\m_\a(\o)={1\/2}\sum_{j,k=1}^d\,M_{jk}\,\o_j
\o_k,
\end{aligned}
\]
where $\S^{d}$ is the $d$-dimensional sphere. Here the linear terms
vanish, since $\l_\a(\vt)$ has an extremum at the point $\vt_0$. The
matrix $M=\{M_{jk}\}_{j,k=1}^d$ is given by
\[\lb{Mjk}
M_{jk}={\pa^2\l_\a(\vt_0)\/\pa\vt_j\pa\vt_k}\,,
\]
and the matrix $m=M^{-1}$ represents a tensor, which is called
\emph{the effective mass tensor} \cite{Ki95}. The effective mass
approximation \er{lam} is a standard approach in solid state
physics. Roughly speaking, in this approach, a complicated
Hamiltonian is replaced by the model Hamiltonian $-{\D\/2m}$\,,
where $\D$ is the Laplacian and $m$ is the so-called effective mass.
We call the quadratic form $\mu_\a(\o)$ \emph{the effective form}.

If a magnetic field is absent, then upper bounds on the effective
masses associated with the ends of each spectral band in terms of
geometric parameters of the graphs were obtained in \cite{KS16}.
Moreover, in the case of the bottom of the spectrum two-sided
estimates on the effective mass in terms of geometric parameters of
the graphs were determined. Now we estimate the effective forms
$\m_\a(\o)$ associated with the ends of each spectral band for the
magnetic Laplacian $\D_\a$.

\begin{theorem}\lb{T2}
Let a band function $\l_\a(\vt)$, $\vt\in\T^d$, have a minimum
(maximum) at some point $\vt_0$ and let $\l_\a(\vt_0)$ be a simple
eigenvalue of $\D_\a(\vt_0)$. Then the effective form $\mu_\a(\o)$ from
\er{lam} satisfies
\[
\lb{em111} \big|\mu_\a(\o)\big|\leq {T_1^2\/\r_\a}+T_2\qqq \forall \
\o\in\S^{d-1},
\]
\[
\text{where}\qqq T_s={1\/s}\,\max_{u\in
V_*}\sum_{\be=(u,v)\in\cS\cup\ul\cS} \|\t
(\be)\|^s\,,\qq s=1,2,
\]
and $\r_\a=\r_\a(\vt_0)$ is the distance between $\l_\a(\vt_0)$ and
the set $\s\big(\D_\a(\vt_0)\big)\setminus\big\{\l_\a(\vt_0)\big\}$, $\t(\be)$ is the index of the edge $\be$ defined by \er{in}, \er{inf},
$\cS$ and $\ul\cS$ are
given by \er{ulS}.
\end{theorem}

\no \textbf{Remarks.}
1) This theorem gives only an upper bound on the effective
form  $\mu_\a(\o)$.  We know low bounds only for the case $\a=0$
\cite{KS16}.

2) Shterenberg \cite{S04}, \cite{S06} considered periodic magnetic Schr\"odinger operators on $\R^d$ and proved that the effective mass tensor can be degenerate for specific magnetic fields, i.e., the matrix $M$ defined by \er{Mjk} is not invertible. In the case of effective masses for magnetic Laplacians on graphs this is an open problem.

\

The paper is organized as follows. In Section \ref{Sec3}
we prove Theorem \ref{TFR2} and Corollary \ref{TCo0} about the decomposition of magnetic Schr\"odinger operators into a constant fiber direct integral with a precise representation of fiber operators. In Section \ref{Sec4} we
prove Theorems \ref{T1}, \ref{Temf} about spectral estimates for magnetic Schr\"odinger operators and Theorem \ref{T2} about estimates on the effective masses of the magnetic Laplacians.
In Section \ref{Sec5} we describe some simple properties of fiber magnetic Laplacians and Schr\"odinger operators and show that the spectral estimates obtained in Theorem \ref{T1} become identities for a specific graph. In the proof we use an example from \cite{KS14}. In Section \ref{Sec5} we also recall some well-known  properties of matrices needed to prove our main results. In Section \ref{Sec6} we consider a more general class of magnetic Laplace and Schr\"odinger operators and briefly formulate similar results for these generalized operators. In this section we also give a factorization of the generalized fiber magnetic Laplacians. This factorization may be crucial for investigation of the bottom of the spectrum of the magnetic Laplacians, for example, for two-sided estimates on the effective mass as it happened in the non-magnetic case (see \cite{KS16}). This section can be read independently on the rest of the paper.

\section{\lb{Sec3} Direct integrals for magnetic Schr\"odinger operators}
\setcounter{equation}{0}

In this section we prove Theorem \ref{TFR2} and Corollary \ref{TCo0}.

\subsection{Floquet decomposition of Schr\"odinger operators}
Recall that we introduce the Hilbert space $\mH$ by \er{Hisp}.
We identify the vertices of the fundamental graph $\G_*=(V_*,\cE_*)$ with the vertices of the periodic graph $\G=(V,\cE)$ from the fundamental vertex set $V_0$.

\begin{theorem}\label{TFD1}
For each 1-form $\a\in\mF_1$ the magnetic Schr\"odinger operator $H_\a=\D_\a+Q$ on $\ell^2(V)$ has the following decomposition into a constant fiber direct integral
\[
\lb{raz}
\begin{aligned}
& \ell^2(V)={1\/(2\pi)^d}\int^\oplus_{\T^d}\ell^2(V_*)\,d\vt ,\qqq
UH_\a U^{-1}={1\/(2\pi)^d}\int^\oplus_{\T^d}\wh H_\a(\vt)d\vt,
\end{aligned}
\]
for the unitary operator $U:\ell^2(V)\to\mH$ defined by
\[
\lb{5001}
(Uf)(\vt,v)=\sum\limits_{m\in\mathbb{Z}^d}e^{-i\lan m,\vt\ran }
f(v+m), \qqq
(\vt,v)\in \T^d\ts V_*,
\]
where $v+m$ denotes the action of $m\in\Z^d$ on $v\in V_*$ and
$\lan\cdot\,,\cdot\ran$ denotes the standard inner product in
$\R^d$.  Here the fiber magnetic Schr\"odinger  operator $\wh
H_\a(\vt)$ and the fiber magnetic Laplacian $\wh \D_\a(\vt)$ are given by
\[\label{Hvt'}
\begin{array}{c}
\wh H_\a(\vt)=\wh \D_\a(\vt)+Q,\\[6pt]
\end{array}
\]
\[
\label{l2.15''}
 \big(\wh \D_\a(\vt)f\big)(v)=\sum_{\be=(v,\,u)\in\cA_*} \big(f(v)-e^{i(\a(\be)+\lan\t(\be),\,\vt\ran)}f(u)\big), \qqq
 v\in V_*,
\]
where $\t(\be)\in\Z^d$ is the edge index defined by \er{in}, \er{inf}.
\end{theorem}

\no {\bf Proof.} Denote by $\ell_{fin}^2(V)$ the
set of all finitely supported functions $f\in \ell^2(V)$. Standard
arguments (see pp. 290--291 in \cite{RS78}) give that $U$ is well
defined on $\ell_{fin}^2(V)$ and has a unique extension to a unitary
operator. For  $f\in \ell_{fin}^2(V)$ the sum \er{5001} is finite
and using the identity $V=\big\{v+m: (v,m)\in V_*\ts\Z^d\big\}$  we
have
$$
\begin{aligned}
&\|Uf\|^2_{\mH}=\int_{\T^d}\|(Uf)(\vt,\cdot)\|_{V_*}^2{d\vt\/(2\pi)^d}\\
&=
 \int_{\T^d}\sum_{v\in V_*}\bigg(\sum\limits_{m\in\Z^d}e^{-i\lan m,\vt\ran }
 f(v+m)\bigg)
 \bigg(\sum\limits_{m'\in\Z^d}e^{i\lan m',\vt\ran }\ol{f(v+m')}\bigg)\,
{d\vt \/(2\pi)^d}\\
\end{aligned}
$$
$$
\begin{aligned}
&=\sum_{v\in V_*}\sum_{m,m'\in\Z^d}
 f(v+m)\ol{f(v+m')}\,\int_{\T^d}e^{-i\lan m-m',\vt\ran}
 {d\vt \/(2\pi)^d}\\
 &=
 \sum_{(v,\,m)\in V_*\ts\Z^d}\big|f(v+m)\big|^2
 =\sum_{v\in V}|f(v)|^2=\|f\|_V^2.
\end{aligned}
$$
Thus, $U$ is  well defined on $\ell_{fin}^2(V)$ and has a unique
isometric extension. In order to prove  that $U$ is onto $\mH$ we
compute $U^*$.   Let $g=\big(g(\cdot,v)\big)_{v\in V_*}\in\mH$,
where $g(\cdot,v) :\T^d\to \C$. We define
\begin{equation}\label{rep}
(U^*g)(v)=\int_{\T^d}e^{i\lan m,\vt\ran }g(\vt,v_*){d\vt \/(2\pi)^d}\,,\qquad
 v=v_*+m\in V,
\end{equation}
where $(v_*,m)\in V_*\ts\Z^d$ is uniquely defined.  A direct computation
gives  that it is indeed the formula for the adjoint of $U$.
Moreover, Parseval's identity for the Fourier series gives
$$
\begin{aligned}
\|U^*g\|^2_V=\sum_{v\in V}\big|(U^*g)(v)\big|^2=\sum_{(v,\,m)\in
V_*\ts\Z^d}\big|(U^*g)(v+m)\big|^2\\= \sum_{(v,\,m)\in
V_*\ts\Z^d}\bigg|\int_{\T^d}e^{i\lan m,\vt\ran }g(\vt,v){d\vt
\/(2\pi)^d}\bigg|^2
\\
=\sum_{v\in V_*}\int_{\T^d}\big|g(\vt,v)\big|^2{d\vt \/(2\pi)^d}\,=
\int_{\T^d}\sum_{v\in V_*}\big|g(\vt,v)\big|^2{d\vt
\/(2\pi)^d}\,=\|g\|_\mH^2.
\end{aligned}
$$
Further, for $f\in \ell_{fin}^2(V)$ and $v\in V_*$ we obtain
$$
\begin{aligned}
(U\D_\a f)(\vt,v)=\sum_{m\in\Z^d}e^{-i\lan m,\vt\ran}(\D_\a f)(v+m)\\
=\sum_{m\in\Z^d}e^{-i\lan m,\vt\ran}
\sum_{\be=(v+m,\,u)\in\cA} \big(f(v+m)-e^{i\a(\be)}f(u)\big)\\
\end{aligned}
$$
$$
\begin{aligned}
=\sum_{\be=(v,\,u)\in\cA_*}\sum_{m\in\Z^d}e^{-i\lan
m,\vt\ran}f(v+m)- \sum_{m\in\Z^d}e^{-i\lan
m,\vt\ran} \sum_{\be=(v,\,u)\in\cA_*}e^{i\a(\be)}f(u+\t(\be)+m)\\
=\sum_{\be=(v,\,u)\in\cA_*}(Uf)(\vt,v)
-\sum_{\be=(v,\,u)\in\cA_*}e^{i(\a(\be)+\lan\t(\be),\vt\ran)}
\sum_{m\in\Z^d}e^{-i\lan m+\t(\be),\vt\ran }f(u+\t(\be)+m)
\end{aligned}
$$
$$
\begin{aligned}
=\sum_{\be=(v,\,u)\in\cA_*}\Big[(Uf)(\vt,v)
-e^{i(\a(\be)+\lan\t(\be),\vt\ran)} (Uf)(\vt,u)\Big]
=\big(\wh\D_\a(\vt)(Uf)(\vt,\cdot)\big)(v).
\end{aligned}
$$
This and the following identity
$$
\begin{aligned}
\big(UQf\big)(\vt,v)=\sum_{m\in\Z^d}e^{-i\lan m,\vt\ran}\big(Q
f\big)(v+m)\\=\sum_{m\in\Z^d}e^{-i\lan m,\vt\ran
}Q(v)f(v+m)=Q(v)(Uf)(\vt,v)
\end{aligned}
$$
yield
$$
(U\D_\a f)(\vt,\cdot)=\wh\D_\a(\vt)(Uf)(\vt,\cdot), \qqq
\big(UQf\big)(\vt,\cdot)=Q\,(Uf)(\vt,\cdot).
$$
Thus, we obtain
$$
UH_\a U^{-1}=U(\D_\a+Q)U^{-1}=
\int_{\T^d}^{\os}\big(\wh\D_\a(\vt)+Q\big)\,{d\vt\/(2\pi)^d}=
\int_{\T^d}^{\os}\wh H_\a(\vt)\,{d\vt\/(2\pi)^d}\,,
$$
which completes the proof. \quad $\BBox$

\subsection{Magnetic fluxes representation}
For a given magnetic field the magnetic potential $\a$ is defined
up to a gauge transformation.
Therefore, in Theorem \ref{TFR2'} we will give a more convenient representation of the fiber magnetic Laplacian $\wh\D_\a(\vt)$ in terms of \emph{magnetic fluxes}.

We fix a vertex $v_0\in V_*$. For each $\vt\in\T^d$ we define  the
function $W:\ell^2(V_*)\ra\R^\n$ as follows: for any vertex
$v\in V_*$, take an oriented path $p=(\be_1,\be_2,\ldots,\be_n)$ on $\G_*$
starting at $v_0$ and ending at $v$ and set
\[\lb{Wvt}
W(v)=\sum_{s=1}^n\big(\a(\be_s)-\a_*(\be_s)\big),
\]
where $\a_*\in\cF_1$ is defined by \er{cat}.

\begin{proposition}
\lb{Tdesp} The value $W(v)$ does not depend on the choice of a path
from $v_0$ to $v$.
\end{proposition}
\no {\bf Proof.} Let $p$ and $q$ be some oriented pathes from
$v_0$ to $v$. We consider the cycle $\cC=p\ul q$, where $\ul
q$ is the inverse path of $q$. Then we have
\[\lb{com1}
\sum_{\be\in\cC}\big(\a(\be)-\a_*(\be)\big)=\sum_{\be\in
p}\big(\a(\be)-\a_*(\be)\big)- \sum_{\be\in
q}\big(\a(\be)-\a_*(\be)\big).
\]
The definition of $\a_*$ gives that for each basic cycle
$\cC_e$ we
have $\sum_{\be\in\cC_e}\a(\be)=\sum_{\be\in\cC_e}\a_*(\be)$, and,
consequently, for each cycle $\cC$ we get $\sum_{\be\in\cC}\a(\be)=\sum_{\be\in\cC}\a_*(\be). $ Combining the
last identity and \er{com1}, we obtain
$$
\sum_{\be\in p}\big(\a(\be)-\a_*(\be)\big)= \sum_{\be\in
q}\big(\a(\be)-\a_*(\be)\big),
$$
which implies that $W(v)$ does not depend on the choice of a path.
\qq $\BBox$

\begin{theorem}
\label{TFR2'} For each $\vt\in\T^d$ the fiber
magnetic Laplacian $\wh\D_\a(\vt)$ defined by \er{l2.15''} is
unitarily equivalent, by a gauge transformation $\cU$ acting in
$\ell^2(V_*)$ and given by
\[\lb{Uvt}
(\cU\,g)(v)=e^{iW(v)}g(v), \qqq g\in\ell^2(V_*),\qqq v\in V_*,
\]
where $W$ is defined by \er{Wvt}, to the operator $\D_\a(\vt)$
given by
\[
\label{gl2.13}
\begin{aligned}
\big(\D_\a(\vt)f\big)(v)=\sum_{\be=(v,\,u)\in\cA_*}  \big(f(v)-e^{i\,\Phi(\be,\vt)}f(u)\big), \qqq v\in V_*,
\end{aligned}
\]
where
\[\lb{phase}
\Phi(\be,\vt)=\a_*(\be)+\lan\t(\be),\vt\ran;
\]
the modified 1-form $\a_*\in\mF_1$ is given by \er{cat},
$\t(\be)$ is the index of the edge $\be$ defined by \er{in}, \er{inf}.
\end{theorem}
\no{\bf Proof.} From \er{Wvt} it follows that
$$
W(u)=W(v)+\a(\be)-\a_*(\be), \qqq \forall\,\be=(v,u)\in\cA_*.
$$
Using this, \er{l2.15''} and \er{Uvt}, we have
\[
\begin{aligned}
\big(\wh\D_\a(\vt)f\big)(v)=\sum_{\be=(v,u)\in\cA_*}
\big(f(v)-e^{i(\a(\be)+\lan\t (\be),\,\vt\ran)}f(u)\big)\\=
\sum_{\be=(v,u)\in\cA_*}
 \big(f(v)-e^{-iW(v)}e^{i\,\Phi(\be,\vt)}e^{iW(u)}f(u)\big)\\
=e^{-iW(v)}\sum_{\be=(v,u)\in\cA_*}
\big(e^{iW(v)}f(v)-e^{i\,\Phi(\be,\vt)}e^{iW(u)}f(u)\big)=
\big(\cU^{-1}\D_{\a}(\vt)\cU f\big)(v),
\end{aligned}
\]
which yields the required statement. \qq $\BBox$

\

\no{\bf Remark.} From Theorem \ref{TFR2'} it follows that each
fiber magnetic Schr\"odinger operator $\wh
H_\a(\vt)=\wh\D_\a(\vt)+Q$, $\vt\in\T^d$, defined by \er{Hvt'} is
unitarily equivalent to the operator $H_\a(\vt)=\D_\a(\vt)+Q$.

\

\no {\bf Proof of Theorem \ref{TFR2}.} This theorem follows from Theorems \ref{TFD1} and \ref{TFR2'}. \qq $\BBox$

\begin{corollary}\label{TCo2}
The $\n\ts\n$ matrix $\D_\a(\vt)=\{\D_{\a,uv}(\vt)\}_{u,v\in V_*}$ associated to the fiber operator $\D_\a(\vt)$ for the magnetic combinatorial Laplacian $\D_\a$ defined by \er{DLO} in the standard orthonormal basis
is given by
\[\lb{Devt}
\D_{\a,uv}(\vt)=\left\{
\begin{array}{ll}
\displaystyle\vk_v-\sum\limits_{\be=(u,u)\in\cA_*}\cos\Phi(\be,\vt), & \textrm{ if } \, u=v\\[20pt]
   \displaystyle-\sum\limits_{\be=(u,v)\in\cA_*}
  e^{-i\,\Phi(\be,\vt)},  & \textrm{ if } \, u\sim v,\qq u\neq v\\[20pt]
 \hspace{20mm} 0, \qq & \textrm{ otherwise}\\
\end{array}\right..
\]
Here  $\n$ is the number of the fundamental graph vertices, ${\vk}_v$ is the degree of the vertex $v$,
$\Phi(\be,\vt)$ is defined by \er{phase}, \er{cat}.
\end{corollary}

\no {\bf Proof.} Let $(\gf_u)_{u\in V_*}$ be the standard orthonormal basis of $\ell^2(V_*)$. Substituting the formula \er{gl2.13} in the identity
$$
\D_{\a,uv}(\vt)=\lan \gf_u,\D_\a(\vt)\gf_v\ran_{V_*}
$$
and using the fact that
for each loop $\be=(u,u)\in\cA_*$ with the phase $\Phi(\be,\vt)$ there exists a loop $\ul\be=(u,u)\in\cA_*$ with the phase $-\Phi(\be,\vt)$ and the identity
$$
e^{-i\,\Phi(\be,\vt)}+e^{i\,\Phi(\be,\vt)}=2\cos \Phi(\be,\vt),
$$
we obtain \er{Devt}. \qq $\BBox$

\

\no {\bf Proof of Corollary \ref{TCo0}.} Due to the connectivity of the $\Z^d$-periodic graph $\G$, on the fundamental graph $\G_*$ there exist $d$ edges $\be_1,\ldots,\be_d$ with linearly independent indices $\t(\be_1),\ldots,\t(\be_d)\in\Z^d$. Then there exists $\vt_0\in\T^d$ satisfying the system of the linear equations
\[\lb{seth}
\a_*(\be_s)+\lan\t(\be_s),\vt_0\ran=0,\qqq s=1,\ldots,d.
\]
If we make the change of variables  $\wt\vt=\vt-\vt_0$, then, using \er{seth} and \er{cat}, for each $\be\in\cA_*$ we have
\begin{multline*}
\a_*(\be)+\lan\t(\be),\vt\ran=\a_*(\be)+\lan\t(\be),\wt\vt+\vt_0\ran \\
=\left\{
\begin{array}{cl}
0,  & \textrm{ if } \, \be\notin(\cS\cup\ul\cS)\\[6pt]
\lan\t(\be),\wt\vt\,\ran,  & \textrm{ if } \, \be\in\{\be_1,\ldots,\be_d,\ul\be_1,\ldots,\ul\be_d\}\\[6pt]
\a_*(\be)+\lan\t(\be),\vt_0\ran+\lan\t(\be),\wt\vt\,\ran, \qq & \textrm{ otherwise}\\
\end{array}\right..
\end{multline*}
Thus,
$$
\a_*(\be)+\lan\t(\be),\vt\ran=\wt\a(\be)+\lan\t(\be),\wt\vt\,\ran,\qqq \forall\,\be\in\cA_*,
$$
where $\wt\a$ is defined by \er{wtcat}, and we obtain \er{wtl2.13}.

Now let $\b=d$. Then
$$
\cS=\{\be_1,\ldots,\be_d\}, \qqq \wt\cS=\cS\sm\{\be_1,\ldots,\be_d\}=\varnothing,\qqq \wt\a=0, \qqq \D_\a(\vt)=\D_0(\wt\vt\,).
$$
This yields that $H_\a$ is unitarily equivalent to $H_0$. \qq $\BBox$

\section{\lb{Sec4} Proof of the main results}
\setcounter{equation}{0}
In this section we prove Theorems \ref{T1} -- \ref{T2}.

\

\no {\bf Proof of Theorem \ref{T1}.} i) We need the following representation of the
Floquet matrix $H_\a(\vt)$, $\vt\in\T^d$:
\[
\label{eq.1}
H_\a(\cdot)=H_\a^0+V_\a(\cdot),\qqq
H_\a^0={1\/(2\pi)^d}\int_{\T^d}H_\a(\vt)d\vt.
\]
From \er{eq.1}, \er{Devt}, \er{inf} and \er{iesp} we deduce that the matrix
$V_\a(\vt)=\{V_{\a,uv}(\vt)\}_{u,v\in V_*}$ has the form
\[
\label{tl2.15} V_{\a,uv}(\vt)=-
\sum\limits_{\be=(u,v)\in\cS\cup\ul\cS}e^{-i\,\Phi(\be,\vt)}.
\]
We define the diagonal matrix $B_\a(\vt)$ by
\[
\lb{eq.2'}
B_\a(\vt)=\diag(B_{\a,u}(\vt))_{u\in V_*},\qqq
B_{\a,u}(\vt)=\sum_{v\in V_*}\big|V_{\a,uv}(\vt)\big|,\qqq \vt\in\T^d.
\]
From \er{tl2.15} we deduce that
\[
\lb{wvv1} \big|V_{\a,uv}(\vt)\big|\leq\big|V_{\a,uv}(0)\big|=\b_{uv}, \qqq
\forall \ (u,v,\vt)\in V_*^2\ts\T^d,
\]
where
\[\lb{buv}
\b_{uv}=\#\{\be\in\cS\cup\ul\cS\mid \be=(u,v)\},
\]
$\#A$ is the number of elements of the set $A$.
Then \er{wvv1} gives
\[
\lb{B0} B_\a(\vt)\leq B_\a(0),\qqq \forall\,\vt\in\T^d.
\]
Then the estimate \er{B0} and Proposition \ref{PrMa}.iii yield
\[
\lb{eq.2} -B_\a(0) \leq -B_\a(\vt)\le V_\a(\vt)\le B_\a(\vt)\leq B_\a(0),\qqq
\forall \vt\in\T^d.
\]
We use some arguments from \cite{Ku15}. Combining
\er{eq.1} and \er{eq.2}, we obtain
$$
H_\a^0-B_\a(0)\le H_\a(\vt)\le H_\a^0+B_\a(0).
$$
Thus, the standard perturbation theory (see Proposition \ref{PrMa}.i) gives
$$
\l_n(H_\a^0-B_\a(0))\leq\l_{\a,n}^-\le \l_{\a,n}(\vt)\leq\l_{\a,n}^+\le \l_n(H_\a^0+B_\a(0)),
\qq \forall \ (n,\vt) \in \N_\n\ts\T^d,
$$
which implies
\[
\lb{009}
 \big|\s(H_\a)\big|\le\sum_{n=1}^{\nu}(\l_{\a,n}^+-\l_{\a,n}^-)\leq
 \sum_{n=1}^\nu\big(\l_n(H_\a^0+B_\a(0))-\l_n(H_\a^0-B_\a(0))\big)=
 2\Tr B_\a(0).
\]
In order to determine $2\Tr B_\a(0)$ we use the relations \er{wvv1}, \er{buv}
and we obtain
\[
\lb{pro} 2\Tr B_\a(0)=2\sum_{u\in V_*}
B_{\a,u}(0)=2\sum_{u,v\in V_*}|V_{\a,uv}(0)|= 2\sum_{u,v\in V_*}\b_{uv}=
2\#(\cS\cup\ul\cS)=4\b.
\]
The estimate \er{eq.7'} follows from \er{009} and \er{pro}.

Now we will prove \er{GEga}.
Since $\l_{\a,1}^-$ and $\l_{\a,\n}^+$ are the lower and upper endpoints of the spectrum $\s(H_\a)$, respectively, using the estimate \er{eq.7'}, we obtain
\[\lb{GEga1}
\sum_{n=1}^s|\g_n(H_\a)|=
\l^+_{\a,\n}-\l_{\a,1}^--\big|\s(H_\a)\big|\geq\l^+_{\a,\n}-\l_{\a,1}^--4\b.
\]
We rewrite the sequence $(Q(v))_{v\in V_*}$ in nondecreasing order
\[
\lb{wtqn}
q_1^\bullet\le q_2^\bullet \le\ldots \le q_\n^\bullet \qq {\rm and \ let } \qq q_1^\bullet=0.
 \]
Here $q_1^\bullet=Q(v_1), q_2^\bullet=Q(v_2),\ldots,q_\n^\bullet=Q(v_\n)$
for some distinct vertices $v_1,v_2,\ldots,v_\n\in V_*$ and without loss of generality we may assume that $q_1^\bullet=0$.

Then Proposition \ref{PrMa}.ii
gives that the eigenvalues of the Floquet matrix $H_\a(\vt)$ for $H_\a=\D_\a+Q$ satisfy
\[
\lb{qq0}
\begin{array}{l}
  q_n^\bullet+\l_{\a,1}^{0-}\le q_n^\bullet+\l_{\a,1}^0(\vt)\le\l_{\a,n}(\vt)\le q_n^\bullet+\l_{\a,\n}^0(\vt)\le q_n^\bullet+\l_{\a,\n}^{0+}\,, \\[6pt]
  \l_{\a,n}^0(\vt)\le\l_{\a,n}(\vt)\le \l_{\a,n}^0(\vt)+q_\n^\bullet\,, \qqq \forall\, (\vt,n)\in\T^d\ts\N_\n.
\end{array}
\]
The first inequalities in \er{qq0} give
\[\lb{qq}
\l_{\a,\n}^+\geq q_\n^\bullet+\l_{\a,1}^{0-}\,,\qqq \l_{\a,1}^-\leq \l_{\a,\n}^{0+}\,,
\]
and, using the second inequalities in \er{qq0}, we have
\[\lb{qq1}
\l_{\a,\n}^{0+}=\max_{\vt\in\T^d}\l_{\a,\n}^{0}(\vt)=
\l_{\a,\n}^{0}(\vt_+)\leq\l_{\a,\n}(\vt_+)\leq \l_{\a,\n}^+\,,
\]
\[\lb{qq2}
\l_{\a,1}^{0-}=\min_{\vt\in\T^d}\l_{\a,1}^{0}(\vt)=\l_{\a,1}^{0}(\vt_-)
\geq\l_{\a,1}(\vt_-)-q_\n^\bullet
\geq\l_{\a,1}^--q_\n^\bullet
\]
for some $\vt_-,\vt_+\in\T^d$.
From \er{qq} -- \er{qq2} it follows that
$$
\l^+_{\a,\n}-\l_{\a,1}^-\geq q_\n^\bullet+\l_{\a,1}^{0-}-\l_{\a,\n}^{0+},\qqq
\l^+_{\a,\n}-\l_{\a,1}^-\geq\l_{\a,\n}^{0+}-\l_{\a,1}^{0-}-q_\n^\bullet,
$$
which yields \er{GEga}.

ii) This item will be proved in Proposition \ref{TG1}.v.
\qq $\BBox$

\

\no {\bf Proof of Theorem \ref{Temf}.} We use the magnetic fluxes
representation given by Theorem \ref{TFR2}. Define the operator
$V_\a(\vt)$, $\vt\in\T^d$, acting on $\C^\n$ by
$$
\begin{aligned}
\D_\a(\vt)=\D_0(\vt)+V_\a(\vt),\qqq  H_\a(\vt)=H_0(\vt)+V_\a(\vt).
\end{aligned}
$$
Here  $\D_0(\vt)$ is the fiber Laplacian and $V_\a(\vt)$ is the fiber magnetic perturbation
operator with the matrix $V_\a(\vt)=\{V_{\a,uv}(\vt)\}_{u,v\in V_*}$ given by
\[
V_{\a,uv}(\vt)=\sum\limits_{\be=(u,v)\in\cS\cup\ul\cS}
e^{-i\lan\t(\be),\,\vt\ran}\big(1-e^{-i\phi_\a(\be)}\big).
\]
Let $\l_{\a,1}(\vt)\le\l_{\a,2}(\vt)\le\ldots\le \l_{\a,\n }(\vt)$ be the eigenvalues of $H_\a(\vt)$.
We have
\[
\begin{aligned}
H_{\a^o}(\vt)=H_\a(\vt)+X_{\a^o,\a}(\vt),\qqq X(\vt)\equiv X_{\a^o,\a}(\vt)=
V_{\a^o}(\vt)-V_{\a}(\vt),
\end{aligned}
\]
where the matrix $X(\vt)=\{X_{uv}(\vt)\}_{u,v\in V_*}$ is given by
 \[\lb{Xuv}
X_{uv}(\vt)=\sum\limits_{\be=(u,v)\in\cS\cup \ul\cS}
  e^{-i\lan\t(\be),\,\vt\ran }\big(e^{-i\phi_\a(\be)}-e^{-i\phi_{\a^o}(\be)}\big).
\]
Then Proposition \ref{PrMa}.ii gives that for each $n\in\N_\n$ we have
\[
\lb{es1} \l_{\a,n}(\vt)+\L_1\leq\l_{\a^o,n}(\vt)\leq\l_{\a,n}(\vt)+\L_\n,
\]
where $\L_1$, $\L_\n$ are defined by \er{emf3'}. From
this we deduce that
$\L_1\le \l_{\a^o,n}^\pm-\l_{\a,n}^{\pm}\le\L_\n$ and
\[
\lb{ees3} |\s_{\a^o,n}|=\l_{\a^o,n}^+-\l_{\a^o,n}^-\le
(\l_{\a,n}^{+}-\l_{\a,n}^{-})
+(\L_\n-\L_1)=|\s_{\a,n}|+(\L_\n-\L_1).
\]
Similar arguments give
\[
\lb{ees4}
\begin{aligned}
|\s_{\a^o,n}|=\l_{\a^o,n}^+-\l_{\a^o,n}^-\geq
(\l_{\a,n}^{+}-\l_{\a,n}^{-})
-(\L_\n-\L_1)=|\s_{\a,n}|-(\L_\n-\L_1).
\end{aligned}
\]
Combining \er{ees3} and \er{ees4}, we obtain \er{emf2}.

We estimate $\L_1$ and $\L_\n$. The standard estimate yields
\[
\|X(\vt)\|\le
\max_{u\in V_*} \sum_{v\in V_*}
\big|X_{uv}(\vt)\big|,\qqq
\max\{|\L_1|,|\L_\n|\}\le \max_{\vt\in\T^d}\|X(\vt)\|.
\]
Using \er{Xuv}, we obtain
\[
|X_{uv}(\vt)|\le \sum\limits_{\be=(u,v)\in\cS\cup\ul\cS}
2\,|\sin x_{\be}|,\qqq
x_{\be}={1\/2}\big(\f_{\a^o}(\be)-\f_{\a}(\be)\big), \qqq \forall\,\vt\in\T^d.
\]
Then we deduce that
\[
\begin{aligned}
\lb{es3.1} \max\{|\L_1|,|\L_\n|\}\le\max_{\vt\in \T^d}\|X(\vt)\|\le
\max_{u\in V_*} \sum_{v\in V_*}\sum_{\be=(u,v)\in\cS\cup\ul
\cS}2\,|\sin x_{\be}|=C_{\a^o,\a},
\end{aligned}
\]
where $C_{\a^o,\a}$ is defined in \er{es3}. The second identity in \er{es30} is a simple consequence of the first identity.
$\BBox$

\

\no \textbf{Proof of Theorem \ref{T2}.}
Let $\p_\a(\vt_0,\cdot)\in\C^\n$ be the normalized eigenfunction,
corresponding to the simple eigenvalue $\l_\a(\vt_0)$. Then the
eigenvalue $\l_\a(\vt)$ and the corresponding normalized
eigenfunction $\p_\a(\vt,\cdot)$ have asymptotics as
$\vt=\vt_0+\ve\o$, $\o\in\S^{d-1}$, $\ve\to0$:
\[
\label{lamp}
\begin{aligned}
&\l_\a(\vt)=\l_\a(\vt_0)+\ve^2\m_\a(\o)+O(\ve^3), \qqq
\p_\a(\vt,\cdot)=\p_{\a,0}+\ve\p_{\a,1}+\ve^2\p_{\a,2}+O(\ve^3),\\
&\textstyle\m_\a(\o)={1\/2}\,\ddot\l_\a(\vt_0+\ve\o)\big|_{\ve=0},\qqq
\p_{\a,0}=\p_\a(\vt_0,\cdot),\\
&\p_{\a,1}=\p_{\a,1}(\o,\cdot)=\dot\p_\a(\vt_0+\ve\o,\cdot)\big|_{\ve=0},\qqq
\textstyle\p_{\a,2}=\p_{\a,2}(\o,\cdot)={1\/2}\,
\ddot\p_\a(\vt_0+\ve\o,\cdot)\big|_{\ve=0},
\end{aligned}
\]
where $\dot u=\pa u/\pa \ve$ and $\S^{d}$ is the $d$-dimensional
sphere.
The Floquet matrix $\D_\a(\vt)$, $\vt\in\T^d$, defined by
\er{Devt} can be
represented in the following form:
\[
\label{geq.1}
\D_\a(\vt)-\l_\a(\vt_0)\1_\n=\D_{\a,0}+\ve\D_{\a,1}(\o)+\ve^2\D_{\a,2}(\o)+O(\ve^3),
\]
as $\vt=\vt_0+\ve\o$, $\ve\rightarrow0$, $\o\in
\S^{d-1}$, where
\[\lb{DDD}
\begin{aligned}
\D_{\a,0}=\D_\a(\vt_0)-\l_\a(\vt_0)\1_\n,\qq  \D_{\a,1}(\o)=\dot
\D_\a(\vt_0+\ve\o)\big|_{\ve=0},\qq \D_{\a,2}(\o)=\textstyle{1\/2}\,\ddot
\D_\a(\vt_0+\ve\o)\big|_{\ve=0},
\end{aligned}
\]
$\1_\n$ is the identity $\n\ts\n$ matrix. The equation
$\D_\a(\vt)\p_\a(\vt,\cdot)=\l_\a(\vt)\p_\a(\vt,\cdot)$ after substitution
\er{lamp}, \er{geq.1} takes the form
\[\lb{gaas}
\begin{aligned}
\big(\D_{\a,0}+\ve\D_{\a,1}(\o)+\ve^2\D_{\a,2}(\o)+O(\ve^3)\big)
\big(\p_{\a,0}+\ve\p_{\a,1}+\ve^2\p_{\a,2}+O(\ve^3)\big)\\=
\big(\ve^2\m_\a(\o)+O(\ve^3)\big)
\big(\p_{\a,0}+\ve\p_{\a,1}+\ve^2\p_{\a,2}+O(\ve^3)\big),
\end{aligned}
\]
where $\p_{\a,0},\p_{\a,1},\p_{\a,2}$ are defined in \er{lamp}. This asymptotics
gives two identities for any $\o\in \S^{d-1}$:
\[\lb{gve1}
\D_{\a,1}(\o)\p_{\a,0}+\D_{\a,0}\p_{\a,1}=0,
\]
\[\lb{gve2}
\D_{\a,2}(\o)\p_{\a,0}+\D_{\a,1}(\o)\p_{\a,1}+\D_{\a,0}\p_{\a,2}=\m_\a(\o)\,\p_{\a,0}.
\]

Using
\er{Devt} we obtain that the entries of the matrices
$\D_{\a,s}(\o)=\{\D^{(\a,s)}_{uv}(\o)\}_{u,v\in V_*}$, $s=1,2$,
defined by \er{DDD}, have the form
\[
\lb{Del1}
\D_{uv}^{(\a,1)}(\o)=i\sum\limits_{\be=(u,\,v)\in\cS\cup\ul\cS}
  \lan\t
(\be),\o\ran\,e^{-i\,\Phi(\be,\vt_0)},
\]
\[
\lb{Del2} \D_{uv}^{(\a,2)}(\o)={1\/2}
\sum\limits_{\be=(u,\,v)\in\cS\cup\ul\cS}\lan\t (\be),\o\ran^2\,
e^{-i\,\Phi(\be,\vt_0)},
\]
for any $\o\in \S^{d-1}$, where $\Phi(\be,\vt)$ is defined by \er{phase}, \er{cat}.

We recall a simple fact that $\D_{\a,1}(\o)\p_{\a,0}$ and $\p_{\a,0}$ are
orthogonal. Indeed, multiplying both sides of \er{gve1} by $\p_{\a,0}$
and using that $\D_{\a,0}\p_{\a,0}=0$, we have $\lan\D_{\a,1}\p_{\a,0},\p_{\a,0}\ran=0$,
which yields $\D_{\a,1}(\o)\p_{\a,0}\perp\p_{\a,0}$.

Let $P_\a$ be the orthogonal projection onto the subspace of $\ell^2(V_*)$ orthogonal to $\p_{\a,0}$. From \er{gve1} we obtain
\[\lb{p_1}
\p_{\a,1}=-(P_\a\D_{\a,0})^{-1}P_\a\D_{\a,1}(\o)\p_{\a,0}.
\]
Multiplying both sides of \er{gve2} by $\p_{\a,0}$, substituting \er{p_1}
and using that $\D_{\a,0}\p_{\a,0}=0$, we have
\[\lb{sec}
\m_\a(\o)=\lan\D_{\a,2}(\o)\p_{\a,0},\p_{\a,0}\ran-
\lan(P_\a\D_{\a,0})^{-1}P_\a\D_{\a,1}(\o)\p_{\a,0},\D_{\a,1}(\o)\p_{\a,0}\ran.
\]
This yields
\[\lb{sec1}
\begin{aligned}
\big|\m_\a(\o)\big|\le\big|\lan\D_{\a,2}(\o)\p_{\a,0},\p_{\a,0}\ran\big|+
\big|\lan(P_\a\D_{\a,0})^{-1}P_\a\D_{\a,1}(\o)\p_{\a,0},\D_{\a,1}(\o)\p_{\a,0}
\ran\big|\\
\leq
\|\D_{\a,2}(\o)\|+\|(P_\a\D_{\a,0})^{-1}P_\a\|\cdot\|\D_{\a,1}(\o)\|^2\leq
\|\D_{\a,2}(\o)\|+{1\/\r_\a}\,\|\D_{\a,1}(\o)\|^2,
\end{aligned}
\]
where $\r_\a=\r_\a(\vt_0)$ is the distance between $\l_\a(\vt_0)$ and $\s\big(\D_\a(\vt_0)\big)\sm\big\{\l_\a(\vt_0)\big\}$. Due to \er{Del1}, \er{Del2}, we have
\[\lb{NN1}
\|\D_{\a,1}(\o)\|\leq\max_{u\in V_*}\sum_{\be=(u,v)\in\cS\cup\ul\cS}
\,\big|\lan\t(\be),\,\o\ran\big|\leq\max_{u\in
V_*}\sum_{\be=(u,v)\in\cS\cup\ul\cS} \|\t(\be)\|=T_1,
\]
\[\lb{NN2}
\|\D_{\a,2}(\o)\|\leq\max_{u\in V_*}\sum_{\be=(u,v)\in\cS\cup\ul\cS}
{\lan\t(\be),\,\o\ran^2\/2} \leq\max_{u\in
V_\ast}\sum_{\be=(u,v)\in\cS\cup\ul\cS}{\|\t (\be)\|^2\/2}=T_2.
\]
Substituting \er{NN1}, \er{NN2} into \er{sec1}, we obtain \er{em111}. \qq
\BBox

\section{Properties of fiber operators and an example \lb{Sec5}}
\setcounter{equation}{0}

In this section we show that the spectral estimates obtained in Theorem \ref{T1} become identities for a specific graph.

\subsection{Properties of fiber operators} We describe some simple properties of fiber magnetic Laplacians and Schr\"odinger operators.

\begin{proposition}\label{TpFo}
For a given 1-form $\a\in \mF_1$, we define another 1-form by   $\wh\a(\be)=-\a(\be)$ for every $\be\in\cA_*$. Then for each $\vt \in\T^d$  the spectra of the fiber magnetic Schr\"odinger operators $\wh H_\a(\vt)$ and $\wh H_{\wh\a}(\vt)$ defined by \er{Hvt'}, \er{l2.15''} satisfy $\s\big(\wh H_{\wh\a}(-\vt)\big)=\s\big(\wh H_{\a}(\vt)\big)$ and, consequently, $\s(H_{\wh\a})=\s(H_{\a})$.
\end{proposition}
\no {\bf Proof.} It is obvious by setting a unitary map $U:\ell^2(V_*)\ra\ell^2(V_*)$ as $U(f)=\ol f$ and using \er{inin}. \qq
$\BBox$

\

A graph is called \emph{bipartite} if its vertex set is divided into two disjoint sets (called \emph{parts} of the graph) such that
each edge connects vertices from distinct parts. A graph is called \emph{regular of degree $\vk_+$} if each its vertex $v$ has
the degree $\vk_v=\vk_+$.

\begin{proposition}\label{TpFo1}
Assume that $\G$ is a periodic regular graph of degree $\vk_+$. Then
the fiber Laplacians $\wh\D_\a(\vt)$ defined by \er{l2.15''} (and, due to unitary equivalence, also the fiber Laplacians $\D_\a(\vt)$ defined by \er{l2.13}) have the following properties.

i) Let we fix any orientation on $\cE_*$ and let $\a$ be a given 1-form. Let  1-form $\wh\a$ be defined as follows:
\[
\wh\a(\be)=\pi-\a(\be) \textrm{ for every } \be\in\cE_*; \qqq \wh\a(\be)=-\wh\a(\ul\be\,) \textrm{ for every }\be\in\cA_*\sm\cE_*.
\]
Then the spectra $\s\big(\wh\D_{\a}(\vt)\big)$ and $\s\big(\wh\D_{\wh\a}(-\vt)\big)$ of the fiber magnetic Laplacians $\wh\D_\a(\vt)$ and $\wh\D_{\wh\a}(-\vt)$ defined by \er{l2.15''}  are symmetric with respect to $\vk_+$ for each $\vt\in\T^d$, that is,
\[\lb{prs1}
\s\big(\vk_+\1-\wh\D_{\a}(\vt)\big)=-\s\big(\vk_+\1-\wh\D_{\wh\a}(-\vt)\big),
\]
where $\1$ is the identity operator. Consequently, $\s(\D_{\a})$ and $\s(\D_{\wh\a})$ are symmetric with respect to $\vk_+$, that is,  $\s(\vk_+\1-\D_{\a})=-\s(\vk_+\1-\D_{\wh\a})$.

ii) Suppose that a fundamental graph $\G_*$ of the graph $\G$ is bipartite. Then  the spectrum $\s\big(\wh\D_{\a}(\vt)\big)$ of the fiber magnetic Laplacian $\wh\D_\a(\vt)$ is symmetric with respect to $\vk_+$ for each $\vt \in\T^d$, that is,
\[
\l\in\s\big(\vk_+\1-\wh\D_{\a}(\vt)\big)\qq \Leftrightarrow \qq -\l\in\s\big(\vk_+\1-\wh\D_{\a}(\vt)\big).
\]
Consequently, $\s(\D_\a)$ is symmetric with respect to $\vk_+$.

\end{proposition}

\no {\bf Proof.}
i) Since the identity $e^{i(\wh\a(\be)-\lan\t
(\be),\,\vt\ran)}=-e^{-i(\a(\be)+\lan\t
(\be),\,\vt\ran)}$ holds true for each $\be\in\cA_*$, it follows that
\[
\begin{aligned}
\big(\big(\vk_+\1-\wh\D_{\wh\a}(-\vt)\big)f\big)(v)=
\sum_{\be=(v,u)\in\cA_*}e^{i(\wh\a(\be)-\lan\t
(\be),\,\vt\ran)}f(u)\\=
-\sum_{\be=(v,u)\in\cA_*}e^{-i(\a(\be)+\lan\t
(\be),\,\vt\ran)}f(u)
=-\overline{\big(\big(\vk_+\1-\wh\D_{\a}(\vt)\big)\ol f\,\big)(v)}.
\end{aligned}
\]
This yields \er{prs1}.

ii) Let $\G_*$ be a bipartite fundamental graph with parts
$V_1$ and $V_2$. We define the unitary operator $U$ on $f\in\ell^2(V_*)$ by
$$
(U f)(v)=\left\{\begin{array}{rl}
f(v),  \ &  {\rm if}\  \ v\in V_1 \\[6pt]
 -f(v), \ &  {\rm if}\  \ v\in V_2  \\
\end{array}\right..
$$
Then we obtain
$$
\big(U^{-1}\big(\vk_+\1-\wh\D_\a(\vt)\big)Uf\big)(v)=
-\big(\big(\vk_+\1-\wh\D_\a(\vt)\big)f\big)(v),
$$
which yields that $\s\big(\wh\D_\a(\vt)\big)$ is symmetric with respect to $\vk_+$. \qq $\BBox$

\

\no \textbf{Remark.} The properties i) and ii) for the magnetic Laplacians on a locally finite graph were proved in \cite{HS01}.

\subsection{Maximal abelian covering.}
We consider a specific class of periodic graphs having some particular properties. Let $\G$ be a $\Z^d$-periodic graph with a fundamental graph $\G_*=(V_*,\cE_*)$ such that $d=\b$, where $\b=\#\cE_*-\#V_*+1$ is the Betti number and $\#A$ is the number of elements of the set $A$. In literature such a periodic graph $\G$ is called  the \emph{maximal abelian covering graph} of the finite graph $\G_*=(V_*,\cE_*)$.  Examples of such graphs include the $d$-dimensional
lattice, the hexagonal lattice.

\textbf{Example.} As an example of the maximal abelian covering
graph we consider a periodic graph, shown in Fig.\ref{ff.10}\emph{a},
and describe the spectrum of the magnetic Laplace and Schr\"odinger operators.

It is known that  $\l_*$ is an eigenvalue
of the Schr\"odinger operator $H_0$ iff $\l_*$ is an eigenvalue of $H_0(\vt)$ for any $\vt\in\T^d$
(see Proposition 4.2 in \cite{HN09}).
Thus, we can define the \emph{multiplicity} of a flat band in the following way:
a flat band $\l_*$  of $H_0$ has multiplicity $m$ iff
$\l_*=\const$ is an eigenvalue of $H_0(\vt)$ of multiplicity $m$ for almost all $\vt\in\T^d$.

\setlength{\unitlength}{1.0mm}
\begin{figure}[h]
\centering
\unitlength 1mm 
\linethickness{0.4pt}
\ifx\plotpoint\undefined\newsavebox{\plotpoint}\fi 
\begin{picture}(120,50)(0,0)

\put(10,10){\line(1,0){40.00}}
\put(10,30){\line(1,0){40.00}}
\put(10,50){\line(1,0){40.00}}
\put(10,10){\line(0,1){40.00}}
\put(30,10){\line(0,1){40.00}}
\put(50,10){\line(0,1){40.00}}

\put(10,10){\circle{1}}
\put(30,10){\circle{1}}
\put(50,10){\circle{1}}

\put(10,30){\circle{1}}
\put(30,30){\circle{1}}
\put(50,30){\circle{1}}

\put(10,50){\circle{1}}
\put(30,50){\circle{1}}
\put(50,50){\circle{1}}

\put(10,10){\vector(1,0){20.00}}
\put(10,10){\vector(0,1){20.00}}

\put(7.5,7.5){$\scriptstyle v_\n$}
\put(1,29){$\scriptstyle v_\n+a_2$}
\put(29,7){$\scriptstyle v_\n+a_1$}
\put(20,8){$\scriptstyle a_1$}
\put(6,20){$\scriptstyle a_2$}
\put(14,26.5){$\scriptstyle v_1$}
\put(23.8,17){$\scriptstyle v_{\nu-2}$}

\put(24.0,12){$\scriptstyle v_{\nu-1}$}
\put(20,26.5){$\scriptstyle v_2$}

\put(10,10){\line(2,3){10.00}}
\put(10,10){\line(3,2){15.00}}
\put(10,10){\line(1,3){5.00}}
\put(10,10){\line(3,1){15.00}}
\put(20.5,21){$\ddots$}
\put(15,25){\circle{1}}
\put(20,25){\circle{1}}
\put(25,15){\circle{1}}
\put(25,20){\circle{1}}

\put(30,10){\line(2,3){10.00}}
\put(30,10){\line(3,2){15.00}}
\put(30,10){\line(1,3){5.00}}
\put(30,10){\line(3,1){15.00}}
\put(40.5,21){$\ddots$}
\put(35,25){\circle{1}}
\put(40,25){\circle{1}}
\put(45,15){\circle{1}}
\put(45,20){\circle{1}}

\put(10,30){\line(2,3){10.00}}
\put(10,30){\line(3,2){15.00}}
\put(10,30){\line(1,3){5.00}}
\put(10,30){\line(3,1){15.00}}
\put(20.5,41){$\ddots$}
\put(15,45){\circle{1}}
\put(20,45){\circle{1}}
\put(25,35){\circle{1}}
\put(25,40){\circle{1}}

\put(30,30){\line(2,3){10.00}}
\put(30,30){\line(3,2){15.00}}
\put(30,30){\line(1,3){5.00}}
\put(30,30){\line(3,1){15.00}}
\put(40.5,41){$\ddots$}
\put(35,45){\circle{1}}
\put(40,45){\circle{1}}
\put(45,35){\circle{1}}
\put(45,40){\circle{1}}
\put(-4,10){(\emph{a})}
\put(64,25){(\emph{b})}

\multiput(80,45)(4,0){5}{\line(1,0){2}}
\multiput(100,25)(0,4){5}{\line(0,1){2}}
\put(80,25){\line(1,0){20.00}}
\put(80,25){\line(0,1){20.00}}
\put(80,25){\circle{1}}
\put(100,25){\circle{1}}

\put(80,45){\circle{1}}
\put(100,45){\circle{1}}

\put(80,25){\vector(1,0){20.00}}
\put(80,25){\vector(0,1){20.00}}

\put(77.0,23.0){$\scriptstyle v_\n$}
\put(76.5,45){$\scriptstyle v_\n$}
\put(101,23){$\scriptstyle v_\n$}
\put(101,45){$\scriptstyle v_\n$}
\put(90,23){$\scriptstyle a_1$}
\put(76,35){$\scriptstyle a_2$}

\put(84,41.5){$\scriptstyle v_1$}
\put(93.8,32.5){$\scriptstyle v_{\nu-2}$}

\put(94.0,27){$\scriptstyle v_{\nu-1}$}
\put(90,41.5){$\scriptstyle v_2$}

\put(80,25){\line(2,3){10.00}}
\put(80,25){\line(3,2){15.00}}
\put(80,25){\line(1,3){5.00}}
\put(80,25){\line(3,1){15.00}}
\put(90.5,36){$\ddots$}
\put(85,40){\circle{1}}
\put(90,40){\circle{1}}
\put(95,30){\circle{1}}
\put(95,35){\circle{1}}
\put(64,10){(\emph{c})}

\put(80,10){\line(1,0){50.00}}
\put(80,9){\line(0,1){2.00}}
\put(95,9){\line(0,1){2.00}}
\put(84,9){\line(0,1){2.00}}
\put(130,9){\line(0,1){2.00}}
\put(85,10){\circle*{1}}

\put(80,9.8){\line(1,0){4.00}}
\put(80,10.2){\line(1,0){4.00}}

\put(95,9.8){\line(1,0){35.00}}
\put(95,10.2){\line(1,0){35.00}}

\put(84.3,11.5){$\scriptstyle1$}
\put(79,6){$\scriptstyle0$}
\put(95,6){$\scriptstyle3$}
\put(123.0,5.5){$\scriptstyle\frac{11+\sqrt{89}}2$}
\put(82.0,5.5){$\scriptstyle\frac{11-\sqrt{89}}2$}
\end{picture}
\vspace{-0.5cm} \caption{\footnotesize  \emph{a}) $\Z^2$-periodic graph $\G$, the vectors $a_1,a_2$ produce an action of $\Z^2$;\quad
\emph{b}) the fundamental graph $\G_*$;\quad \emph{c}) the spectrum of the
Laplacian ($\n=3$).} \label{ff.10}
\end{figure}
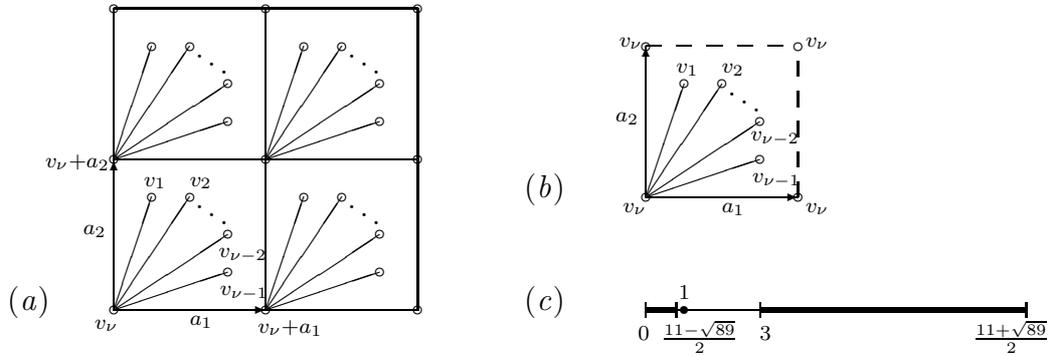

\begin{proposition}\lb{TG1}
Let $\dL^d_*$ be the fundamental graph of the
$d$-dimension lattice $\dL^d$ and let $v_\n$ be the unique vertex of $\dL^d_*$.
Let $\G_*$ be obtained from $\dL^d_*$ by adding $\n-1\geq1$ vertices
$v_1,\ldots,v_{\n-1}$ and $\n-1$ unoriented edges
$(v_1,v_\n),\ldots,(v_{\n-1},v_\n)$ (see
Fig.\ref{ff.10}b). Let $\a\in \mF_1$ be a given 1-form. Then
$\G$ is a maximal abelian covering graph of $\G_*$ and satisfies

i) The spectrum of the magnetic Laplacian $\D_\a$ on the periodic graph $\G$ has the form
\[
\lb{acs1}
\s(\D_\a)=\s(\D_0)=\s_{ac}(\D)\cup\s_{fb}(\D),\qqq \s_{fb}(\D)=\s_2(\D)=\{1\},
\]
where the flat band $\s_2(\D)=\{1\}$ has multiplicity $\n-2$ and  $\s_{ac}(\D)$
has only two bands $\s_1(\D)$ and $\s_3(\D)$ given by
\[
\lb{acs2}
\begin{array}{c}
\s_{ac}(\D)=\s_1(\D)\cup\s_3(\D), \\ [6 pt]
\s_1(\D)=\big[\,0,\textstyle x-\sqrt{x^2-4d}\;\big],\qq
\s_3(\D)=\big[\,\n,\textstyle x+\sqrt{x^2-4d}\;\big],\qqq
x={\n+4d\/2}\,.
\end{array}
\]

ii) The spectrum of the magnetic Schr\"odinger operator $H_\a=\D_\a+Q$ on $\G$ has the form
\[
\lb{acs333}
\s(H_\a)=\s(H_0)=\bigcup_{n=1}^\n\big[\l_n(0),\l_n(\vt_\pi)\big], \qqq \vt_\pi=(\pi,\ldots,\pi)\in\T^d.
\]

iii) Let $q_\n=Q(v_\n)=0$ and let all other values of  the potential $q_1=Q(v_1),\ldots,q_{\n-1}=Q(v_{\n-1})$ at the
vertices of the fundamental graph $\G_*$ be distinct. Then
$\s(H_0)=\s_{ac}(H_0)$, i.e., $\s_{fb}(H_0)=\varnothing$.

iv) Let among the numbers $q_1,\ldots,q_{\n-1}$ there exist a value
$q_*$ of multiplicity $m$. Then the spectrum of the Schr\"odinger
operator $H_0$ on $\G$ has the flat band $q_*+1$ of multiplicity
$m-1$.

v) The Lebesgue measure of the spectrum of the magnetic Schr\"odinger operators $H_\a$ on $\G$ satisfies
\[
\lb{sp2}
\textstyle|\s(H_\a)|=4d
\]
and the estimates \er{eq.7'} and the first estimate in \er{GEga} become identities.
\end{proposition}

\no {\bf Proof.} The fundamental graph $\G_*$ consists of $\n$  vertices
$v_1,v_2,\ldots,v_\n$; ${\n-1}$ unoriented edges
$(v_1,v_\n),\ldots,(v_{\n-1},v_\n)$ and $d$ unoriented loops in the
vertex $v_\n$. Since
\[\lb{beqd}
\b=\#\cE_*-\#V_*+1=(\n-1+d)-\n+1=d,
\]
the graph $\G$ is a maximal
abelian covering graph of $\G_*$ and, due to Corollary \ref{TCo0},
$\s(H_\a)=\s(H_0)$.

Items i) -- v) for the case $\a=0$ were proved in \cite{KS14} (Proposition 7.2).
Combining \er{sp2} and \er{beqd} we obtain $|\s(H_\a)|=4\b$, i.e., the estimates \er{eq.7'} become identities. \qq
\BBox

\subsection{Well-known  properties of matrices}
\lb{App}
We recall some well-known  properties of matrices (see e.g.,
\cite{HJ85}). Denote by $\l_1(A)\leq\ldots\leq\l_\n(A)$ the
eigenvalues of a self-adjoint $\n\ts\n$ matrix $A$, arranged in
increasing order, counting multiplicities.

\begin{proposition}\lb{PrMa}
i) Let $A,B$ be self-adjoint $\nu\ts\nu$ matrices and let
$B\geq0$. Then the eigenvalues $\l_n(A)\leq\l_n(A+B)$ for all
$n\in\N_\n$ (see Corollary 4.3.3 in \cite{HJ85}).

ii) Let $A,B$ be self-adjoint $\nu\ts\nu$ matrices. Then for each
$n\in\N_\n$ we have
$$
\l_n(A)+\l_1(B)\leq\l_n(A+B)\leq\l_n(A)+\l_\n(B)
$$
(see Theorem 4.3.1 in \cite{HJ85}).

iii) Let $V=\{V_{jk}\}$ be a self-adjoint $\n\ts \n$  matrix, for
some $\n<\iy$ and let $ B=\diag \{B_1,\ldots, B_\n\}$,
$B_j=\sum\limits_{k=1}^\n |V_{jk}|$. Then the following estimates
hold true:
\[
\begin{aligned}
\lb{V} &-B\le V\le B
\end{aligned}
\]
(see \cite{K13}).
\end{proposition}

\section{Generalized magnetic Schr\"odinger operators \lb{Sec6}}
\setcounter{equation}{0}

\subsection{Generalized magnetic Laplacians on periodic graphs.}
In this section we deal with a more general class of magnetic Laplacians. These generalized magnetic Laplacians on finite and infinite graphs are considered in \cite{CTT11}, \cite{HS99a}, \cite{HS99b}, \cite{HS01}, \cite{LLPP15}, \cite{S94}. We define two positive weights on $\G$
\[
m_V:V\rightarrow(0,\iy),\qqq m_\cA:\cA\rightarrow(0,\iy)
\]
such that
\[
m_V(v+\mm)=m_V(v), \qqq m_\cA(\be+\mm)=m_\cA(\be)=m_\cA(\ul\be)
\]
for all $(v,\be,\mm)\in V\ts\cA\ts\Z^d$.
We consider the weighted Hilbert space
\[\lb{ellV}
\ell^2(V,m_V)=\Big\{f: V\ra\C \mid \sum_{v\in V}m_V(v)|f(v)|^2<\infty\Big\},
\]
equipped with the inner product
\[\lb{ipV}
\lan f,g\ran_V=\sum_{v\in V}m_V(v)f(v)\ol{g(v)}.
\]
For each 1-form $\a\in\mF_1$ we define the \emph{discrete magnetic Laplace operator} $\D_\a$ on
$f\in\ell^2(V,m_V)$ by
\[
\lb{DLO1}
 \big(\D_\a f\big)(v)={1\/m_V(v)}\sum_{\be=(v,u)\in\cA}m_\cA(\be) \big(f(v)-e^{i\a(\be)}f(u)\big), \qqq
 v\in V.
\]

\no \textbf{Remark.} If $\a=0$, then $\D_0$ is just the discrete Laplacian $\D$:
\[
\big(\D f\big)(v)={1\/m_V(v)}\sum_{\be=(v,u)\in\cA}m_\cA(\be) \big(f(v)-f(u)\big), \qqq v\in V.
\]

It is well known (see \cite{HS99a}, \cite{HS99b}, \cite{HS01}) that \emph{the magnetic Laplacian $\D_\a$ is a bounded self-adjoint operator on $\ell^2(V,m_V)$ and its spectrum $\s(\D_\a)$ is a closed subset in $[0,2\vk_+]$, where $\vk_+$ is defined by
\[\lb{vkpl1}
\vk_+=\sup_{v\in V}{1\/m_V(v)}\sum_{\be=(v,u)\in\cA}m_\cA(\be).
\]
Here the sum is taken over all oriented edges starting at the vertex $v$.}

\

We present typical magnetic Laplacians:

1)  \textbf{The magnetic combinatorial Laplacian.} If we set $m_V(v)=1$
for each vertex $v\in V$ and  $m_\cA(\be)=1$ for each edge
$\be\in\cA$, then the \emph{magnetic combinatorial Laplacian} is
expressed by \er{DLO} and is discussed in Sections \ref{Sec1}--\ref{Sec4}. This
magnetic Laplacian and  corresponding Schr\"odinger operators are
considered in \cite{B13}, \cite{DM06}, \cite{LL93}.

  2) \textbf{The magnetic transition operator.} Let $p:\cA\rightarrow(0,1]$ be a transition probability such that
\[
\sum_{\be=(v,u)\in\cA}p(\be)=1, \qqq \forall \,v\in V.
\]
Moreover, let $p$ be $m_V$-symmetric, that is $m_V(v)p(\be)=m_V(u)p(\ul\be\,)$ for each oriented edge $\be=(v,u)$. If we set $m_\cA(\be)=m_V(v)p(\be)$, then the magnetic Laplace operator is expressed by
\[\lb{dml}
\D_\a=\1-T_{p,\a},\qqq \big(T_{p,\a} f\big)(v)=\sum_{\be=(v,u)\in\cA}p(\be)\,e^{i\a(\be)}f(u),
\]
where $T_{p,\a}$ is the \emph{magnetic transition operator with respect to $p$} and $\1$ is the identity operator. The magnetic Laplacian $\D_\a=\1-T_{p,\a}$ is considered in \cite{HS99a}, \cite{HS99b}, \cite{HS01}.

3) {\bf The magnetic normalized Laplacian.} This Laplacian is obtained from \er{dml} if we set $p(\be)={1\/\vk_v}$ for each $\be=(v,u)\in\cA$. Then $m_V(v)=\vk_v$, $m_\cA(\be)=1$ and the \emph{magnetic normalized Laplacian} is expressed as follows:
\[\lb{mno}
\big(\D_\a f\big)(v)=f(v)-{1\/\vk_v}\sum_{\be=(v,u)\in\cA}e^{i\a(\be)}f(u).
\]

\subsection{Main results for generalized Schr\"odinger operators}
In this subsection we generalize the results formulated in Section \ref{Sec2} (Theorems \ref{TFR2}, \ref{T1} and \ref{T2}) for the magnetic Schr\"odinger operator $H_\a=\D_\a+Q$ with the Laplacian $\D_\a$ defined by \er{DLO1}.

We introduce the Hilbert space
\[
\mH=L^2\Big(\T^d,{d\vt
\/(2\pi)^d}\,,\cH\Big)=\int_{\T^d}^{\os}\cH\,{d\vt \/(2\pi)^d}\,, \qq \cH=\ell^2(V_*,m_{V_*}),
\]
i.e., a constant fiber direct integral equipped with the norm
$$
\|g\|^2_{\mH}=\int_{\T^d}\|g(\vt,\cdot)\|_{V_*}^2\frac{d\vt
}{(2\pi)^d}\,,
$$
where the function $g(\vt,\cdot)\in\cH$ for almost all $\vt\in\T^d$.

\begin{theorem}\label{GTFD1}
For each 1-form $\a\in\mF_1$ the magnetic Schr\"odinger operator $H_\a=\D_\a+Q$ with the Laplacian $\D_\a$ defined by \er{DLO1} on $\ell^2(V,m_V)$ has the following decomposition into a constant fiber direct integral
\[
\lb{Graz}
\begin{aligned}
& \ell^2(V,m_V)={1\/(2\pi)^d}\int^\oplus_{\T^d}\ell^2(V_*,m_{V_*})\,d\vt ,\qqq
\mU H_\a\mU^{-1}={1\/(2\pi)^d}\int^\oplus_{\T^d}H_\a(\vt)d\vt,
\end{aligned}
\]
for some unitary operator $\mU:\ell^2(V,m_V)\to\mH$. Here the fiber magnetic Schr\"odinger  operator $H_\a(\vt)$ and the  fiber  magnetic Laplacian $\D_\a(\vt)$ are given by
\[\label{GHvt'}
H_\a(\vt)=\D_\a(\vt)+Q, \qqq \forall\,\vt\in \T^d,
\]
\[
\label{Ggl2.13}
\big(\D_\a(\vt)f\big)(v)={1\/m_{V_*}(v)}\sum_{\be=(v,\,u)\in\cA_*}m_{\cA_*}(\be)  \big(f(v)-e^{i(\a_*(\be)+\lan\t(\be),\vt\ran)}f(u)\big), \qqq v\in V_*,
\]
where the modified 1-form $\a_*\in\mF_1$ is defined by \er{cat},
$\t(\be)$ is the index of the edge $\be$ defined by \er{in}, \er{inf},
and $\lan\cdot\,,\cdot\ran$ denotes the standard inner product in $\R^d$.
\end{theorem}

\begin{theorem}
\lb{GT1g}
The Lebesgue measure $|\s(H_\a)|$ of the spectrum of the magnetic Schr\"odinger operator $H_\a=\D_\a+Q$ with the Laplacian $\D_\a$ defined by \er{DLO1} satisfies
\[
\lb{Geq.7}
|\s(H_\a)|\le \sum_{n=1}^{\n}|\s_n(H_\a)|\le 2\wh\b,
\]
where
\[\lb{Gbeta}
\wh\b=\sum\limits_{\be=(u,v)\in\cS\cup\ul\cS}
{m_{\cA_*}(\be)\/\big(m_{V_*}(u)m_{V_*}(v)\big)^{1/2}}\,,
\]
$\n$ is the number of the fundamental graph vertices, $\cS$ and $\ul\cS$ are
defined by \er{ulS}.
\end{theorem}

\begin{theorem}\lb{T2.G}
Let a band function $\l_\a(\vt)$, $\vt\in\T^d$, of the magnetic Laplacian
$\D_\a$ defined by \er{DLO1} have a minimum
(maximum) at some point $\vt_0$ and let $\l_\a(\vt_0)$ be a simple
eigenvalue of $\D_\a(\vt_0)$. Then the effective form $\mu_\a(\o)$ from
\er{lam} satisfies
\[
\lb{em11} \big|\mu_\a(\o)\big|\leq {T_1^2\/\r_\a}+T_2\qqq \forall\,
\o\in\S^{d-1},
\]
\[
\text{where}\qqq T_s={1\/s}\,\max_{u\in
V_*}\sum_{\be=(u,v)\in\cS\cup\ul\cS}{m_{\cA_*}(\be)\,\|\t
(\be)\|^s\/\big(m_{V_*}(u)m_{V_*}(v)\big)^{1/2}}\,,\qq s=1,2,
\]
where $\r_\a=\r_\a(\vt_0)$ is the distance between $\l_\a(\vt_0)$ and the set
$\s\big(\D_\a(\vt_0)\big)\setminus\big\{\l_\a(\vt_0)\big\}$, $\t(\be)$ is the index of the edge $\be$ defined by \er{in}, \er{inf}; $\cS$ and $\ul\cS$ are defined by \er{ulS}.
\end{theorem}

The proof of these results are similar to the proof of Theorems \ref{TFR2}, \ref{T1} and \ref{T2}.

\subsection{Factorization of fiber magnetic Laplacians}
We introduce the Hilbert space
\[
\ell^2(\cA_*,m_{\cA_*})=\{\phi: \cA_*\ra\C\mid \phi(\ul\be)=-\phi(\be)
\textrm{ for } \be\in\cA_* \qq \textrm{ and } \qq \lan
\phi,\phi\ran_{\cA_*}<\iy\},
\]
where the inner product is given by
\[\lb{incA}
\lan
\phi_1,\phi_2\ran_{\cA_*}={1\/2}\sum_{\be\in\cA_*}m_{\cA_*}(\be)\phi_1(\be)\ol{\phi_2(\be)}.
\]

For each $\vt\in\T^d$ we define the operator
$\na_\a(\vt):\ell^2(V_*,m_{V_*})\ra \ell^2(\cA_*,m_{\cA_*})$ by
\[\lb{navt}
\begin{aligned}
& \big(\na_\a(\vt)
f\big)(\be)=e^{-i\,\Phi(\be,\vt)/2}f(v)-e^{i\,\Phi(\be,\vt)/2}f(u),\qqq \forall\,f\in \ell^2(V_*,m_{V_*}),\\[6pt]
& \textrm{where}\qq \be=(v,u),\qqq
\Phi(\be,\vt)=\a_*(\be)+\lan\t(\be),\vt\ran;
\end{aligned}
\]
the modified 1-form $\a_*\in\mF_1$ is given by \er{cat},
$\t(\be)$ is the index of the edge $\be$ defined by \er{in}, \er{inf}.

\begin{theorem}\label{FOFa}
i) For each $\vt\in\T^d$ the conjugate operator
$\na^*_\a(\vt):\ell^2(\cA_*,m_{\cA_*})\ra \ell^2(V_*,m_{V_*})$ has
the form
\[\lb{coop}
(\na^*_\a(\vt)
\phi)(v)=\sum_{\be=(v,u)\in\cA_*}{m_{\cA_*}(\be)\/m_{V_*}(v)}\,e^{i\,\Phi(\be,\vt)/2}\phi(\be),\qqq
\forall\,\phi\in \ell^2(\cA_*,m_{\cA_*}).
\]

ii) For each $\vt\in\T^d$ the fiber magnetic Laplacian
$\D_\a(\vt)$ defined by \er{Ggl2.13} satisfies
\[\lb{fact}
\D_\a(\vt)=\na^*_\a(\vt)\na_\a(\vt).
\]

iii) If some $\vt\in\T^d$ satisfies
\[\lb{cqm}
\Phi(\be,\vt)=0 ,\qqq \forall\, \be\in\cS,
\]
where $\cS$ is defined in \er{ulS}, then the rank of the operator $\na_\a(\vt)$ is equal to $\n-1$.
Otherwise, the rank of the operator $\na_\a(\vt)$ is equal to $\n$, where $\n$ is the number of the fundamental graph vertices.

iv) For each $\vt\in\T^d$ the quadratic form $\lan \D_\a(\vt)
f,f\ran_{V_*}$ associated with the fiber magnetic Laplacian
$\D_\a(\vt)$ is given by
\[
\lan \D_\a(\vt)
f,f\ran_{V_*}={1\/2}\sum_{\be=(v,u)\in\cA_*}m_{\cA_*}(\be)\,
\big|f(v)-e^{i\,\Phi(\be,\vt)}f(u)\big|^2.
\]
\end{theorem}
\no {\bf Proof.} Let $\vt\in\T^d$, $f\in\ell^2(V_*,m_{V_*})$,
$\phi\in\ell^2(\cA_*,m_{\cA_*})$. Using \er{mpin} and \er{inin} we
have
\[\lb{bein}
\Phi(\ul\be,\vt)=-\Phi(\be,\vt),\qqq
\forall(\be,\vt)\in\cA_*\ts\T^d.
\]

i) Due to \er{incA}, \er{navt}, \er{bein}, we have
\[\lb{coop1}
\begin{aligned}
\lan\na_\a(\vt)
f,\phi\ran_{\cA_*}={1\/2}\sum_{\be\in\cA_*}m_{\cA_*}(\be)\big(\na_\a(\vt)
f\big)(\be)\ol{\phi(\be)}
\\={1\/2}\sum_{\be=(v,u)\in\cA_*}m_{\cA_*}(\be)
\big(e^{-i\,\Phi(\be,\vt)/2}f(v)-e^{i\,\Phi(\be,\vt)/2}f(u)\big)\ol{\phi(\be)}
\\={1\/2}\sum_{\be=(v,u)\in\cA_*}m_{\cA_*}(\be)\,
e^{-i\,\Phi(\be,\vt)/2}f(v)\ol{\phi(\be)}
-{1\/2}\sum_{\be=(v,u)\in\cA_*}m_{\cA_*}(\be)\,
e^{i\,\Phi(\be,\vt)/2}f(u)\ol{\phi(\be)}
\\={1\/2}\sum_{\be=(v,u)\in\cA_*}m_{\cA_*}(\be)\,
e^{-i\,\Phi(\be,\vt)/2}f(v)\ol{\phi(\be)}
+{1\/2}\sum_{\ul\be=(u,v)\in\cA_*}m_{\cA_*}(\ul\be)\,
e^{-i\,\Phi(\ul\be,\vt)/2}f(u)\ol{\phi(\ul\be)}
\\=\sum_{\be=(v,u)\in\cA_*}m_{\cA_*}(\be)\,
e^{-i\,\Phi(\be,\vt)/2}f(v)\ol{\phi(\be)}.
\end{aligned}
\]
On the other hand, due to \er{ipV}, \er{coop}, we obtain
\[\lb{coop2}
\begin{aligned}
\lan f,\na^*_\a(\vt) \phi\ran_{V_*}=\sum_{v\in
V_*}m_{V_*}(v)f(v)\ol{(\na^*_\a(\vt) \phi)(v)}\\=\sum_{v\in
V_*}m_{V_*}(v)f(v)\sum_{\be=(v,u)\in\cA_*}{m_{\cA_*}(\be)\/m_{V_*}(v)}\,e^{-i\,\Phi(\be,\vt)/2}\,\ol\phi(\be)
=\sum_{\be=(v,u)\in\cA_*}m_{\cA_*}(\be)\,
e^{-i\,\Phi(\be,\vt)/2}f(v)\ol{\phi(\be)}.
\end{aligned}
\]
Comparing \er{coop1} and \er{coop2}, we get the required statement.

ii) Using \er{navt}, \er{coop}, we obtain
\[
\begin{aligned}
\big(\na^*_\a(\vt)\na_\a(\vt)f\big)(v)=
\sum_{\be=(v,u)\in\cA_*}{m_{\cA_*}(\be)\/m_{V_*}(v)}\,e^{i\,\Phi(\be,\vt)/2}
\big(\na_\a(\vt)f\big)(\be)\\
=\sum_{\be=(v,u)\in\cA_*}{m_{\cA_*}(\be)\/m_{V_*}(v)}\,e^{i\,\Phi(\be,\vt)/2}
\big(e^{-i\,\Phi(\be,\vt)/2}f(v)-e^{i\,\Phi(\be,\vt)/2}f(u)\big)\\=
\sum_{\be=(v,u)\in\cA_*}{m_{\cA_*}(\be)\/m_{V_*}(v)}\,
\big(f(v)-e^{i\,\Phi(\be,\vt)}f(u)\big)=\big(\D_\a(\vt)f\big)(v),\qqq
\forall v\in V_*.
\end{aligned}
\]

iii) We omit the proof, since it is similar to the proof of
Proposition 2.4.ii in \cite{KS16}.

iv) From \er{fact}, \er{incA}, \er{navt} it follows that
\[
\begin{aligned}
\lan \D_\a(\vt) f,f\ran_{V_*}=\lan
\na_\a(\vt)f,\na_\a(\vt)f\ran_{\cA_*}=
\|\na_\a(\vt)f\|^2_{\cA_*}\\={1\/2}\sum_{\be=(v,u)\in\cA_*}m_{\cA_*}(\be)\,
\big|\big(\na_\a(\vt)f\big)(\be)\big|^2={1\/2}\sum_{\be=(v,u)\in\cA_*}m_{\cA_*}(\be)\,
\big|e^{-i\,\Phi(\be,\vt)/2}f(v)-e^{i\,\Phi(\be,\vt)/2}f(u)\big|^2\\=
{1\/2}\sum_{\be=(v,u)\in\cA_*}m_{\cA_*}(\be)\,
\big|f(v)-e^{i\,\Phi(\be,\vt)}f(u)\big|^2.
\end{aligned}
\]
$\BBox$

\

\no \textbf{Remarks.} 1) The magnetic Laplacian $\D_\a$ defined by \er{DLO1} has the following factorization (see \cite{HS99a},
\cite{HS99b}, \cite{HS01}):
\[
\D_\a=\na^*_\a\na_\a,
\]
where the operator $\na_\a:\ell^2(V,m_V)\ra \ell^2(\cA,m_\cA)$ is
given by
\[
(\na_\a f)(\be)=e^{-i\a(\be)/2}f(v)-e^{i\a(\be)/2}f(u),\qqq
\forall\,f\in \ell^2(V,m_V),\qq\textrm{where}\qq \be=(v,u).
\]
The conjugate operator $\na^*_\a:\ell^2(\cA,m_\cA)\ra \ell^2(V,m_V)$
has the form
\[
(\na^*_\a\,
\phi)(v)=\sum_{\be=(v,u)\in\cA}{m_\cA(\be)\/m_V(v)}\,e^{i\a(\be)/2}\phi(\be),\qqq
\forall\,\phi\in \ell^2(\cA,m_\cA).
\]
The quadratic form $\lan \D_\a f,f\ran_V$ associated with the
magnetic Laplacian $\D_\a$ is given by
\[
\lan \D_\a f,f\ran_V={1\/2}\sum_{\be=(v,u)\in\cA}m_\cA(\be)\,
\big|f(v)-e^{i\a(\be)}f(u)\big|^2.
\]

2) The quasimomentum $\vt$ satisfying \er{cqm} may or may not exist.
For example, if $\#\cS=d$, then such $\vt\in\T^d$
exists and is unique.


\medskip

\footnotesize
\textbf{Acknowledgments. \lb{Sec7}}  Evgeny Korotyaev was
supported by the RSF grant  No. 15-11-30007. Natalia Saburova was
supported by the RFFI grant No. 16-01-00087.
Finally, we would like to thank referees for thoughtful comments that helped us to improve the manuscript.

\end{document}